%% file: noisygs_techreport.tex
\documentclass[10pt]{article}

\RequirePackage{amsmath}
\RequirePackage{amsfonts}
\RequirePackage{amssymb}
\RequirePackage{booktabs}
\RequirePackage{hyperref}
\RequirePackage{multirow}
\RequirePackage{algorithm}
\RequirePackage{algpseudocode}
\RequirePackage{graphicx}
\RequirePackage{xcolor}
\RequirePackage{amsthm}
\RequirePackage{amsmath}
\RequirePackage{amsfonts}
\RequirePackage{amssymb}
\RequirePackage{hyperref}
\RequirePackage{multirow}
\RequirePackage{algorithm}
\RequirePackage{algpseudocode}
\RequirePackage{graphicx}
\RequirePackage{xcolor}
\RequirePackage{enumitem}
\allowdisplaybreaks
\usepackage{pgfplots}
\pgfplotsset{compat=1.18}

\usepackage[
    backend=biber,       
    style=numeric,       
    sorting=nyt,         
    isbn=true,           
    url=true,            
    doi=true,            
    eprint=true,         
    backref=true         
]{biblatex}

\usepackage{dpr_fec,ifthen,dsfont}

\newtheorem{example}{Example}[section]

\usepackage{isetechreport}
\def\reportyear{26T}
\def\reportno{006}
\def\revisionno{1}
\def\originaldate{April 1, 2026}
\def\revisiondate{\today}
\coraltrue
\cvcrfalse

\begin{document}

\title{\input{title}}

\author{Albert S.~Berahas\thanks{E-mail: albertberahas@gmail.com}}
\affil{Department of Industrial and Operations Engineering, University of Michigan}
\author{Frank E.~Curtis\thanks{E-mail: frank.e.curtis@lehigh.edu}}
\author{Lara Zebiane\thanks{E-mail: lara.zebiane@lehigh.edu}}
\affil{Department of Industrial and Systems Engineering, Lehigh University}
\titlepage

\maketitle

\begin{abstract}
  \input{abstract}
\end{abstract}

\input{body}

\appendix
\input{appendix}

\printbibliography

\end{document}

%% file: title.tex
Gradient Sampling for Noisy Nonsmooth Nonconvex Optimization

%% file: abstract.tex
An algorithm is proposed, analyzed, and tested for minimizing locally Lipschitz objective functions that may be nonconvex and/or nonsmooth.  The algorithm is built upon the gradient-sampling methodology and designed specifically for cases when objective and generalized gradient values are subject to bounded errors that might not be controllable by the algorithm (say, due to a computational budget).  A novel assumption pertaining to such errors is introduced that decouples the noise that is present due to typical derivative estimates from noise that is due to misidentification of segments of differentiability.  It is proved that, with probability one, either the sequence of objective values will decrease without bound or the algorithm will terminate finitely at a point at which a measure of stationarity is below a threshold that is determined by, and vanishes with, the error bounds for the objective and generalized gradients.  The results of several numerical experiments are presented, most notably for a fixed-order $H_\infty$ controller design problem and a PDE-constrained optimal control problem, where in both cases errors arise due to computational noise. The results show that the algorithm optimizes reliably despite these errors.

%% file: body.tex
\newcommand{\ball}[3]{\Bmbb_{\leq #2}^{#1}(#3)}
\newcommand{\gengradient}{\eth}
\newcommand{\fapprox}{\ftilde}
\newcommand{\gapprox}{\gtilde}
\newcommand{\ones}{\mathds{1}}
\newcommand{\epsilonls}{\epsilon_{\rm ls}}

\section{Introduction}\label{sec.introduction}

Many problems in control and systems theory require the minimization of functions that are locally Lipschitz, and might not be either convex or smooth.  Examples include the design of fixed-order controllers, where a worst-case performance measure such as the $H_\infty$-norm of a closed-loop transfer function is minimized \cite{ApkaNoll2006,BurkHenrLewiOver2006b,BurkHenrLewiOver2006a}, and optimal control problems governed by (partial) differential equations in which sparsity-promoting control costs appear in the objective~\cite{CasasHerzogWachsmuth2012,Stadler2009}.  A feature that such problems share is that each evaluation of the objective function, and a (generalized) gradient of it, requires the solution of an inner problem---e.g., a level-set or eigenvalue computation \cite{BoydBalaKaba1989,BruiStei1990}, or state and adjoint equations---that in practice is solved by an iterative method only up to a tolerance.  The values that are available for the optimization are therefore subject to errors that are deterministic and bounded, and cannot be driven to zero at a reasonable cost.  Moreover, near a point of nondifferentiability, such errors can cause the wrong active segment of the objective function (e.g., the wrong frequency, eigenvalue, or grid point) to be identified, in which case the error in a computed generalized gradient is large no matter what tolerance is used for the inner solver.  This paper considers the design and analysis of an optimization algorithm with convergence guarantees in such settings.

In settings in which the objective is differentiable and the gradient is continuous, classical derivative-based algorithms are effective, and their behavior in the presence of bounded errors is by now well understood.  General locally Lipschitz optimization, on the other hand, requires specialized algorithms.  One of the most successful frameworks for solving such challenging problems is the gradient sampling (GS) methodology, which first arose in the context of optimizing the stability of dynamical systems~\cite{BurkLewiOver2002} and has since been employed in software for fixed-order controller design \cite{BurkHenrLewiOver2006b}.  Analyzed in foundational work by Burke, Lewis, and Overton \cite{BurkLewiOver05}, and later refined by Kiwiel \cite{Kiwi07}, this methodology involves evaluating gradients at points that are sampled randomly in a neighborhood around each iterate to construct a search direction that reliably approximates a direction of descent. Subsequent research has significantly expanded this framework to include adaptive sampling strategies, quasi-Newton acceleration, and sequential quadratic programming methods for solving constrained problems; see, e.g., the articles \cite{CurtLi22,CurtMitcOver17,CurtOver12,CurtQue13,CurtQue15,CurtRobiZhou20}, as well as the survey~\cite{BurkCurtLewiOverSimo20}.

Despite these efforts in algorithm design and analysis, previously proposed GS methods have, with few exceptions, relied on exact objective function and generalized gradient values being available.  The exceptions of which we are aware are the derivative-free variants in \cite{HareNuti2013,Kiwi2010}, which employ exact objective values along with gradient approximations whose accuracy is controlled by the algorithm, and the trust-region variant for convex problems in \cite{ShanSu2026}. However, for many modern settings, an algorithm only has access to inexact values due to numerical approximation, simulation-based modeling, or stochastic estimation procedures.  In the context of a GS method for general locally Lipschitz settings, the case of errors that are merely bounded, and not controllable by the algorithm, has not been explored (see~\S\ref{sec.literature}), even though errors are particularly disruptive in such settings as they corrupt the sampled gradient information, which leads to unreliability in the search direction computation. This gap in the literature motivates the design of a GS method that handles errors in objective function and generalized gradient values while maintaining rigorous convergence guarantees that are consistent with the noiseless setting.

\subsection{Contributions}

We propose and analyze a GS-based algorithm for solving nonconvex and nonsmooth optimization problems that remains robust in the presence of errors in both the objective function and generalized gradient values when these errors are only assumed to be bounded (with known bounds), and so may be biased and uncontrollable by the algorithm. Consistent with previously proposed gradient-based methods for noisy smooth optimization \cite{BeraByrdNoce2019,BeraCaoSche2021}, our approach perturbs the sufficient decrease condition in the backtracking line search to prevent the failure of the search in the presence of noise. We prove that, with probability one, the algorithm either generates objective function values that decrease without bound or terminates finitely at an iterate satisfying a stationarity condition whose accuracy is tied directly to the magnitudes of the evaluation errors.  Our guarantee has two notable features.  First, the stationarity threshold depends on the square root of the bound on the objective errors, which matches the dependence in state-of-the-art guarantees for noisy smooth optimization.  Second, although the guarantee is stated under an assumption about the probability with which good sample sets are generated, we show in an appendix that this assumption can be replaced by a lemma when the generalized gradient approximations are continuous over an open set of full measure, which is the natural counterpart of the differentiability assumption that underlies the analysis of GS methods in the noiseless setting.  \textit{Critical to our work is a novel assumption about the noise in the generalized gradient information that decouples the noise that may be due to shifts in the domain of the optimization variables from the noise that may be due to typical errors in the derivative evaluations; see Assumption~\ref{ass.g} on page~\pageref{ass.g} and our subsequent discussions of it.} To the best of our knowledge, the proposed method is the first GS-based method for solving nonconvex problems that explicitly incorporates bounded, not-necessarily controllable evaluation errors in both objective and generalized gradient information while preserving convergence guarantees.

\subsection{Literature Review}\label{sec.literature}

One major distinction in the literature on noisy optimization is whether the noise is deterministic, meaning that two calls to an evaluation oracle with the same input always produce the same value, or stochastic.  The former may arise, e.g., as computational noise \cite{BeraByrdNoce2019,OztoByrdNoce23} generated through a numerical simulation, whereas the latter may arise, e.g., in stochastic-gradient-type methods for solving machine learning and other data-driven problems \cite{BottCurtNoce18}.

There is a long history of stochastic/noisy subgradient-type methods for solving unconstrained nonsmooth minimization problems; see, e.g., \cite{BianHachSche2022,BoltLePauw2023,DaviDrusKakaLee2020,ErmoNork1998,NediBert2010,Nork1986,Nork2021,Nurm1974,Rusz2020,Rusz2021,SoloZavr1998} and also \cite{DaviDrus2019}.  Those with guarantees for nonconvex problems require somewhat restrictive assumptions, such as weak convexity, whereas GS methods are applicable to a larger class of locally Lipschitz problems.  The method in~\cite{DaviDrusLeePadm2022}, which builds off of~\cite{ZhanLinJegeSraJadb2020}, is related to a GS approach and, like ours, aims to compute a point that is approximately stationary in the sense of Goldstein~\cite{Gold1977}; however, the guarantees for these methods presume exact generalized gradients or unbiased stochastic gradient estimates, whereas the errors that we consider are deterministic and may be biased. A GS method that employs approximate function evaluations is proposed in \cite{BoskCortMart2026}, although the approach is distinct to our bounded-noise model.

More closely related to our setting are proximal-bundle methods, many of which, like GS methods, have theoretical guarantees for locally Lipschitz minimization under loose assumptions.  For such methods that employ inexact objective and/or generalized gradient values, see \cite{deOlSagaLema2014,HareSagaSolo2016,HertUlbr2019,Hint2001,HuangNiuLinYinYuan2023,Kiwi2006,Noll2013,Solo2003}; in particular, \cite{HareSagaSolo2016,HertUlbr2019} propose and analyze proximal-bundle methods for solving nonconvex and nonsmooth minimization problems when the objective and generalized-gradient values may be subject to noisy evaluations, as we presume in our context.

Inexact evaluations are also a central theme in the literature on optimization problems governed by differential equations, where the accuracy of state and adjoint solves is often adapted by the optimization algorithm \cite{HeinVice2002,KourHeinRidzvanB2014,ZiemUlbr2011}. Such approaches presume that errors can be controlled and that any nonsmooth part of the objective has a structure that can be handled through a proximal mapping \cite{BaraKour2023}.  By contrast, we consider errors that are bounded but not necessarily controllable, and objectives whose nonsmoothness has no exploitable structure; e.g., see our control example in~\S\ref{sec.pde}.

GS methods are arguably simpler to implement than proximal-bundle methods, especially for solving nonconvex problems, since the latter require techniques known as downshifting or tilting that introduce parameters that need to be tuned carefully.  By contrast, our GS-based algorithm involves relatively few parameters, and our theoretical guarantees state explicitly the bounds within which they are to be chosen.

\subsection{Notation}

We use $\N{} := \{1,2,\dots\}$ to denote the set of positive integers.  We use $\R{}$ to denote the set of real numbers, $\R{}_{\geq r}$ (respectively, $\R{}_{>r}$) to denote the set of real numbers greater than or equal to (respectively, greater than) $r \in \R{}$, $\R{n}$ to denote the set of real $n$-vectors, and $\R{m \times n}$ to denote the set of real $m$-by-$n$ matrices. For any positive integer $m \in \N{}$, we define the set $[m] := \{1,\dots,m\}$.

Given a nonempty set $\Xcal \subseteq \R{n}$, we denote the distance function from points in the domain $\R{n}$ to the set $\Xcal$ by $\dist_\Xcal : \R{n} \to \R{}$, which is defined by
\bequationNN
  \dist_{\Xcal}(v) = \inf_{x \in \Xcal} \|x - v\|_2\ \ \text{for all}\ \ v \in \R{n}.
\eequationNN
If $\Xcal \subseteq \R{n}$ is nonempty, closed, and convex, then by the Projection Theorem~\cite[Theorem~6.5]{CurtRobi25} there exists a unique projection of any $x \in \R{n}$ onto $\Xcal$, call it $\proj_{\Xcal}(x)$.  Given $\epsilon \in \R{}_{\geq0}$ and $x \in \R{n}$, we denote the closed Euclidean ball in $\R{n}$ centered at $x$ with radius $\epsilon$ as $\ball{n}{\epsilon}{x} := \{\xbar \in \R{n} : \|\xbar - x\|_2 \leq \epsilon\}$. Given two nonempty sets $\Xcal \subseteq \R{n}$ and $\bar\Xcal \subseteq \R{n}$, we denote $\Xcal + \bar\Xcal = \{x + \xbar \in \R{n} : x \in \Xcal\ \text{and}\ \xbar \in \bar\Xcal\}$.

Suppose that a function $f : \R{n} \to \R{}$ is locally Lipschitz in the sense that, at each point $x \in \R{n}$, there exist $\epsilon_{f,x} \in \R{}_{>0}$ and $L_{f,x} \in \R{}_{>0}$ such that
\bequation\label{eq.locally_Lipschitz}
  |f(\xbar) - f(\xhat)| \leq L_{f,x} \|\xbar - \xhat\|_2\ \ \text{for all}\ \ (\xbar,\xhat) \in \ball{n}{\epsilon_{f,x}}{x} \times \ball{n}{\epsilon_{f,x}}{x}.
\eequation
By Rademacher's theorem \cite{EvansGari2015,Rade1919}, this implies that $f$ is continuous over $\R{n}$ and differentiable almost everywhere in $\R{n}$, i.e., the set in $\R{n}$ over which it is not differentiable has a Lebesgue measure of zero \cite{Clar83,CurtRobi25}.  Following \cite{Clar83}, we denote the Clarke generalized directional derivative function for $f$ as $f^\circ : \R{n} \times \R{n} \to \R{}$ with
\bequationNN
  f^\circ(x,d) = \limsup_{\substack{\xbar \to x \\ \alpha \searrow 0}} \frac{f(\xbar + \alpha d) - f(\xbar)}{\alpha}\ \ \text{for all}\ \ (x,d) \in \R{n} \times \R{n}.
\eequationNN
Subsequently, we denote the Clarke generalized gradient mapping for $f$---a set-valued mapping---by $\gengradient f : \R{n} \rightrightarrows \R{n}$, which is defined for all $x \in \R{n}$ by
\begin{align*}
  \gengradient f(x) :=&\ \{g \in \R{n} : f^\circ(x,d) \geq g^Td\ \text{for all}\ d \in \R{n}\} \\
  =&\ \conv(\{ g \in \R{n} : \{\nabla f(x_k)\} \to g\ \text{for some}\ \{x_k\} \to x\ \text{with}\ \{x_k\} \subset \Dcal_{\nabla f} \}).
\end{align*}
Here, $\Dcal_{\nabla f} \subseteq \R{n}$ is the full-measure set in $\R{n}$ over which $f$ is differentiable.  For the equivalence indicated by the latter equation above, see, e.g., \cite{Clar83,CurtRobi25,RockWets1998}.

\subsection{Organization}

In \S\ref{sec.algorithm}, we state our problem of interest, our assumptions about the evaluation errors, and our proposed algorithm.  We present our convergence analysis in~\S\ref{sec.analysis}, the results of numerical experiments in \S\ref{sec.numerical}, and concluding remarks in~\S\ref{sec.conclusion}.  In Appendix~\ref{app.sampling}, we show that one of the assumptions that is used in our analysis can be replaced by a lemma when the generalized gradient approximations are continuous over an open set of full measure.  Appendices~\ref{app.rosenbrock} and~\ref{app.nn} contain two additional experiments, and Appendix~\ref{app.details} discusses further experimental details.

\section{Problem and Algorithm Statements}\label{sec.algorithm}

Our algorithm is designed to solve (approximately) the problem of minimizing
an objective function $f : \R{n} \to \R{}$, as in
\bequation\label{prob.opt}
  \min_{x \in \R{n}}\ f(x),
\eequation
where $f$, along with a corresponding approximation function $\fapprox$ that is employed by the algorithm, at least satisfy the following assumption (recall \eqref{eq.locally_Lipschitz}).

\bassumption\label{ass.f}
  The objective $f : \R{n} \to \R{}$ is locally Lipschitz, and there exist $\epsilon_f \in \R{}_{\geq0}$ and a (Borel) measurable approximation function $\fapprox : \R{n} \to \R{}$ such that
  \bequationNN
    |\fapprox(x) - f(x)| \leq \epsilon_f\ \ \text{for all}\ \ x \in \R{n}.
  \eequationNN
\eassumption

Under Assumption~\ref{ass.f}, $f$ may be nonconvex and/or nonsmooth. The existence of $\epsilon_f$ and $\fapprox$ such that $|\fapprox(x) - f(x)| \leq \epsilon_f$ for all $x \in \R{n}$ is a common assumption for noisy smooth optimization; see, e.g., \cite{BeraByrdNoce2019,OztoByrdNoce23}. Let us emphasize that any two evaluations of~$\fapprox$ at the same~$x \in \R{n}$ always yield the same value, namely,~$\fapprox(x)$, although $\fapprox$ is not assumed to be Lipschitz continuous; indeed, it may even be discontinuous.

As is generally the case, we assume that at any point $x \in \R{n}$ it would be intractable to compute the entire set of generalized gradients $\gengradient f(x)$, although it would be tractable to compute (or approximate) an element of this set. Along these lines, for our proposed algorithm, we assume that a generalized gradient approximation function is available that, for each $x$, returns a single vector as a generalized gradient approximation. To state our assumption about this function, we introduce two objects. First, let us introduce the $\epsilon$-generalized gradient mapping $\gengradient_\epsilon f : \R{n} \rightrightarrows \R{n}$ as defined by Goldstein~\cite{Gold1977}, which is defined for all $\epsilon \in \R{}_{\geq0}$ and for any $x \in \R{n}$ by
\bequationNN
  \gengradient_\epsilon f(x) := \cl\(\conv\( \Ucal_{f,\epsilon}(x) \)\),\ \ \text{where}\ \ \Ucal_{f,\epsilon}(x) := \bigcup_{\xbar \in \ball{n}{\epsilon}{x}} \gengradient f(\xbar).
\eequationNN
At the same time, given a function $\gapprox : \R{n} \to \R{n}$ (which will represent our generalized gradient approximation function), let us introduce for any $\epsilon \in \R{}_{>0}$ the mapping $\tilde\Gcal_{\epsilon} : \R{n} \rightrightarrows \R{n}$ as being defined for all $x \in \R{n}$ by
\bequationNN
  \tilde\Gcal_{\epsilon}(x) := \cl(\conv(\{\gapprox(\xbar) : \xbar \in \ball{n}{\epsilon}{x}\})).
\eequationNN
We now make the following assumption about the generalized gradient approximation function $\gapprox$ that is available to the algorithm.  (The factor of 1/3 in part (b) of the assumption could be generalized to any factor in (0,1). We merely use a factor of 1/3 as a natural choice that is consistent with the test that is employed in the algorithm's line search to decide whether the search should be abandoned; see upcoming Lemma~\ref{lem.alpha_bound}.)

\bassumption\label{ass.g}
  There exist $\epsilon_x \in \R{}_{\geq0}$, $\epsilon_g \in \R{}_{\geq0}$, and a (Borel) measurable approximation function $\gapprox : \R{n} \to \R{n}$ such that the following hold.
  \benumerate
    \item[(a)] $\gapprox(x) \in \gengradient_{\epsilon_x} f(x) + \ball{n}{\epsilon_g}{0}$ for all $x \in \R{n}$.
    \item[(b)] $\gengradient_{\epsilon/3} f(x) \subseteq \tilde\Gcal_\epsilon(x) + \ball{n}{\epsilon_g}{0}$ for all $x \in \R{n}$ and $\epsilon \in \R{}_{>0}$ with $\epsilon \geq \epsilon_x$.
  \eenumerate
\eassumption

\noindent
Let us emphasize that, under Assumption~\ref{ass.g}, any two evaluations of $\gapprox$ at the same point~$x$ always yield the same value $\gapprox(x)$.  Otherwise, Assumption~\ref{ass.g} requires some explanation as it is a unique contribution to the literature on noisy optimization. The main idea of the assumption, which allows our algorithm and our analysis of it to reduce to the noiseless case when $(\epsilon_f,\epsilon_x,\epsilon_g) \to (0,0,0)$, is to decouple the noise in the generalized gradients that may be present due to certain shifts in the domain from noise that may be present due to more typical derivative evaluation errors.

Assumption~\ref{ass.g} is best understood through a few examples.

\begin{example}
  Consider the graphs pertaining to $f : \R{} \to \R{}$ with $f(x) = |x|$ on the upper-left of Figure~\ref{fig.example} along with the noisy example on the upper-right. In this example, the noise in the objective values is consistent throughout the domain, as is commonly expected in the context of noisy optimization. As for the noise in the generalized gradient values, this particular example turns out to be relatively straightforward since the jump in the graph of $\gapprox$ occurs at $x=0$, i.e., the unique point at which the true set of generalized gradients is not a singleton. Considering any $\epsilon_x \in \R{}_{\geq0}$, one finds that $\gapprox(0) \approx 1 \in [-1,1] + \ball{n}{\epsilon_g}{0} = \gengradient_{\epsilon_x} f(0) + \ball{n}{\epsilon_g}{0}$ for small $\epsilon_g \in \R{}_{>0}$, and for any $x \neq 0$ one finds that $\gapprox(x) \in \{\nabla f(x)\} + \ball{n}{\epsilon_g}{0} \subseteq \gengradient_{\epsilon_x} f(x) + \ball{n}{\epsilon_g}{0}$ for small $\epsilon_g \in \R{}_{>0}$. Thus, the inclusion in part~(a) of Assumption~\ref{ass.g} holds for any $\epsilon_x \in \R{}_{\geq0}$, small $\epsilon_g \in \R{}_{>0}$, and all $x \in \R{}$. At the same time, now for any $\epsilon_x \in \R{}_{\geq0}$ and small $\epsilon_g \in \R{}_{>0}$, one finds for any $\epsilon \in \R{}_{>0}$ with $\epsilon \geq \epsilon_x$ that $\gengradient_{\epsilon/3} f(0) = [-1,1] \subseteq \tilde\Gcal_\epsilon(0) + \ball{n}{\epsilon_g}{0}$ since $\tilde\Gcal_\epsilon(0)$ captures values of $\gapprox$ both to the left and to the right of the origin.  All that remains is to verify that there exist $\epsilon_x \in \R{}_{\geq0}$ and $\epsilon_g \in \R{}_{>0}$ such that the inclusion in part~(b) of the assumption holds for all $x \neq 0$ and $\epsilon \in \R{}_{>0}$ with $\epsilon \geq \epsilon_x$.  Consider any $\epsilon_x \in \R{}_{\geq0}$.  If such $(x,\epsilon)$ yields an interval $[x-\epsilon/3,x+\epsilon/3]$ that does not contain the origin, then $\gengradient_{\epsilon/3} f(x) = \{\nabla f(x)\}$, which is contained in $\{\gapprox(x)\} + \ball{n}{\epsilon_g}{0}$ for small $\epsilon_g \in \R{}_{>0}$, and so contained in $\tilde\Gcal_\epsilon(x) + \ball{n}{\epsilon_g}{0}$ for small $\epsilon_g \in \R{}_{>0}$.  On the other hand, if $[x-\epsilon/3,x+\epsilon/3]$ contains the origin, then the larger interval $[x-\epsilon,x+\epsilon]$ contains the origin in its interior, in which case one finds that $\gengradient_{\epsilon/3} f(x) = [-1,1]$ is a subset of $\tilde\Gcal_{\epsilon} (x) + \ball{n}{\epsilon_g}{0}$ for small $\epsilon_g \in \R{}_{>0}$. Putting everything together, one finds that Assumption~\ref{ass.g} holds with any $\epsilon_x \in \R{}_{\geq0}$ and small $\epsilon_g \in \R{}_{>0}$.
\end{example}

\begin{figure}[ht]
  \centering
  \pgfmathsetseed{1}
  \begin{tikzpicture}
    \begin{axis}[
      width=0.5\textwidth,
      axis lines = middle,
      xlabel = $x$,
      ymin = -1.5, ymax = 2.5,
      xmin = -2.5, xmax = 2.5,
      grid = both,
      grid style = {dashed, gray!30},
      tick label style = {font=\small},
      label style = {font=\small},
      ytick = {0, 2}, 
      extra y ticks = {-1, 1},
      extra y tick labels = {},
      ]
      \node [blue, anchor=south, font=\small] at (axis cs:1.5,1.8) {$f(x)$};
      \addplot [domain=-2:2, samples=100, color=blue, thick] {abs(x)};
      \draw [magenta, thick] (axis cs:-2,-1) -- (axis cs:0,-1);
      \draw [magenta, thick] (axis cs:0,1) -- (axis cs:2,1);
      \draw [magenta, thick] (axis cs:0,-1) -- (axis cs:0,1);
      \addplot [only marks, mark=*, magenta, mark size=1.5pt] coordinates {(0,-1) (0,1)};
      \node [magenta, anchor=south, font=\small] at (axis cs:0.5,1) {$\gengradient f(x)$};
      \node [anchor=east, xshift=-2pt, font=\small] at (axis cs:0,1) {1};
      \node [anchor=west, xshift=2pt, font=\small] at (axis cs:0,-1) {-1};
    \end{axis}
  \end{tikzpicture}
  \hspace{0.8cm}
  \begin{tikzpicture}
    \begin{axis}[
      width=0.5\textwidth, 
      axis lines = middle,
      xlabel = $x$,
      ymin = -1.5, ymax = 2.5,
      xmin = -2.5, xmax = 2.5,
      grid = both,
      grid style = {dashed, gray!30},
      tick label style = {font=\small},
      label style = {font=\small},
      ytick = {0, 2}, 
      extra y ticks = {-1, 1},
      extra y tick labels = {},
      ]
      \node [blue, anchor=south, font=\small] at (axis cs:1.5,1.8) {$\fapprox(x)$};
      \addplot [domain=-2:2, samples=150, color=blue, thick] {abs(x) + 0.12*rand};
      \addplot [domain=-2:0, samples=80, magenta, thick] {-1 + 0.12*rand}; 
      \addplot [domain=0:2, samples=80, magenta, thick] {1 + 0.12*rand};
      \draw [magenta, thick, dashed] (axis cs:0,-1) -- (axis cs:0,1);
      \addplot [only marks, mark=*, magenta, mark size=1.5pt] coordinates {(0,1)};
      \node [magenta, anchor=south, font=\small] at (axis cs:0.7,1.1) {$\gapprox(x)$};
      \node [anchor=east, xshift=-2pt, font=\small] at (axis cs:0,1) {1};
      \node [anchor=west, xshift=2pt, font=\small] at (axis cs:0,-1) {-1};
    \end{axis}
  \end{tikzpicture}
  \vspace{0.4cm}
  
  \begin{tikzpicture}
    \begin{axis}[
      width=0.5\textwidth, 
      axis lines = middle,
      xlabel = $x$,
      ymin = -1.5, ymax = 2.5,
      xmin = -2.5, xmax = 2.5,
      grid = both,
      grid style = {dashed, gray!30},
      tick label style = {font=\small},
      label style = {font=\small},
      ytick = {0, 2}, 
      extra y ticks = {-1, 1},
      extra y tick labels = {},
      extra x ticks = {0.2},
      extra x tick labels = {$\delta$},
      every extra x tick/.style = {tick label style = {xshift=4pt}},
      ]
      \node [blue, anchor=south, font=\small] at (axis cs:1.7,1.8) {$\fapprox(x)$};
      \addplot [domain=-2:2, samples=100, color=blue, thick] {abs(x-0.2)};
      \draw [magenta, thick] (axis cs:-2,-1) -- (axis cs:0.2,-1);
      \draw [magenta, thick] (axis cs:0.2,1) -- (axis cs:2,1);
      \draw [magenta, thick, dashed] (axis cs:0.2,-1) -- (axis cs:0.2,1);
      \addplot [only marks, mark=*, magenta, mark size=1.5pt] coordinates {(0.2,1)};
      \node [magenta, anchor=south, font=\small] at (axis cs:0.7,1) {$\gapprox(x)$};
      \node [anchor=east, xshift=-2pt, font=\small] at (axis cs:0,1) {1};
      \node [anchor=east, xshift=-2pt, font=\small] at (axis cs:0,-1) {-1};
    \end{axis}
  \end{tikzpicture}
  \hspace{0.8cm}
  \begin{tikzpicture}
    \begin{axis}[
      width=0.5\textwidth, 
      axis lines = middle,
      xlabel = $x$,
      ymin = -1.5, ymax = 2.5,
      xmin = -2.5, xmax = 2.5,
      grid = both,
      grid style = {dashed, gray!30},
      tick label style = {font=\small},
      label style = {font=\small},
      ytick = {0, 2}, 
      extra y ticks = {-1, 1},
      extra y tick labels = {},
      extra x ticks = {0.2},
      extra x tick labels = {$\delta$},
      every extra x tick/.style = {tick label style = {xshift=4pt}},
      ]
      \node [blue, anchor=south, font=\small] at (axis cs:1.5,1.8) {$\fapprox(x)$};
      \addplot [domain=-2:2, samples=150, color=blue, thick] {abs(x-0.2) + 0.15*rand};
      \addplot [domain=-2:0.2, samples=80, magenta, thick] {-1 + 0.15*rand}; 
      \addplot [domain=0.2:2, samples=80, magenta, thick] {1 + 0.15*rand};
      \draw [magenta, thick, dashed] (axis cs:0.2,-1) -- (axis cs:0.2,1);
      \addplot [only marks, mark=*, magenta, mark size=1.5pt] coordinates {(0.2,1)};
      \node [magenta, anchor=south, font=\small] at (axis cs:0.7,1.1) {$\gapprox(x)$};
      \node [anchor=east, xshift=-2pt, font=\small] at (axis cs:0,1) {1};
      \node [anchor=west, xshift=2pt, font=\small] at (axis cs:0,-1) {-1};
    \end{axis}
  \end{tikzpicture}
  \caption{Upper-left: absolute value function and its generalized-gradient mapping.  Upper-right: the same subject to bounded errors where the point of nondifferentiability of $f$ is preserved exactly at the jump in the graph of~$\gapprox$.  Here, as is typical in noisy optimization, the mapping~$\gapprox$ does \textit{not} correspond to the derivative function of $\fapprox$; rather, the function $\gapprox$ approximates elements of the generalized-gradient mapping directly. For this particular example on the upper-right, Assumption~\ref{ass.g} holds with arbitrary $\epsilon_x \geq 0$ and positive $\epsilon_g \approx 0$.  Lower-left: a shifted absolute value function $\fapprox(x) = |x - \delta|$ for some $\delta \in \R{}_{>0}$ and its generalized-gradient mapping $\gapprox$. Assumption~\ref{ass.g} holds for this example with $\epsilon_x > 3\delta/2$ and $\epsilon_g = 0$.  Lower-right: the same mappings subject to bounded errors, where errors in the generalized-gradient mapping are due to misidentification of the point of nondifferentiability as well as typical derivative errors.  Assumption~\ref{ass.g} holds with $\epsilon_x > 3\delta/2$ and $\epsilon_g \approx 0$.}
  \label{fig.example}
\end{figure}

The example illustrated in the upper-right of Figure~\ref{fig.example} represents an idealized case in which the jump in $\gapprox$ occurs exactly at the point at which the true generalized-gradient mapping does not yield a singleton. Consider next another idealized situation in which the noise in both the objective and generalized gradient values is due solely to a shift in the formulation of the approximate objective function, and beyond the consequences of this shift the generalized gradient values have no additional error.

\begin{example}
  Consider the example illustrated in the lower-left of Figure~\ref{fig.example}.  The noise in the objective is again consistent throughout the domain, as can be seen in the shifted graph.  As for $\gapprox$, however, there is a large difference, e.g., between $\{\gapprox(x)\} = \{-1\}$ and $\gengradient f(x) = \{\nabla f(x)\} = \{1\}$ for all $x \in (0,\delta)$.  That said, one can deduce that Assumption~\ref{ass.g} holds for appropriately chosen $\epsilon_x \in \R{}_{\geq0}$ and $\epsilon_g = 0$. Let $\epsilon_g = 0$ be fixed for the remainder of the discussion of this example.  With respect to part~(a), one finds for any $\epsilon_x \geq \delta$ that $x < \delta$ gives $\gapprox(x) = -1 \in \gengradient f(\xbar) \subseteq \gengradient_{\epsilon_x} f(x)$ for some $\xbar \leq 0$ with $|x - \xbar| \leq \delta \leq \epsilon_x$, whereas $x \geq \delta$ gives $\gapprox(x) = 1 = \nabla f(x) \in \gengradient_{\epsilon_x} f(x)$.  Hence, the inclusion in part~(a) holds for all $x \in \R{}$.  Now to show that part (b) also holds, let us choose any $\epsilon_x > 3\delta/2$ (say, arbitrarily close to $3\delta/2$) and any $\epsilon \geq \epsilon_x$.  If $x < -\epsilon/3$ or $x > \delta + \epsilon/3$, then $\gengradient_{\epsilon/3} f(x) = \{\nabla f(x)\} = \{\gapprox(x)\} \subseteq \tilde\Gcal_{\epsilon}(x)$.  On the other hand, for any $x \in [-\epsilon/3,\delta + \epsilon/3]$ one has $\delta < 2\epsilon_x/3 \leq 2\epsilon/3$, from which it follows that the distance from any such $x$ to both the origin and $\delta$ is strictly less than $\epsilon$.  Thus,
  \bequationNN
    \tilde\Gcal_\epsilon(x) = [-1,1] \supseteq \gengradient_{\epsilon/3} f(x)\ \text{for all}\ x \in [-\epsilon/3,\delta + \epsilon/3].
  \eequationNN
  Overall, Assumption~\ref{ass.g} holds with any $\epsilon_x > 3\delta/2$ and $\epsilon_g = 0$.  This example demonstrates a critical feature of Assumption~\ref{ass.g}. That is, if errors in the generalized-gradient mapping occur essentially due to misidentification of the points of nondifferentiability, then such errors can be captured by $\epsilon_x > 0$ while the value for $\epsilon_g$ can remain small---to account only for errors in derivative values similarly as in the context of noisy smooth optimization, where a typical assumption is that the approximate gradient mapping $\gapprox$ yields $\|\gapprox(x) - \nabla f(x)\|_2 \leq \epsilon_g$ for all $x \in \R{n}$ \cite{BeraByrdNoce2019,OztoByrdNoce23}.
\end{example}

To complete this discussion, let us provide one last example.

\begin{example}
  Consider the realistic example illustrated on the lower-right of Figure~\ref{fig.example}. Here, noise in the generalized-gradient mapping is present due to both misidentification of the point of nondifferentiability and typical noise in the gradients at points of differentiability.  Combining the insights from the previous examples, one finds that Assumption~\ref{ass.g} holds with $\epsilon_x > 3\delta/2$ and $\epsilon_g \approx 0$. It is worthwhile to emphasize that the assumption can also be seen to hold if the approximation function~$\gapprox$ experienced \textit{multiple} jumps (down to approximately $-1$ and up to approximately~$1$) in an interval containing the origin, as long as $\epsilon_x > 0$ is large enough with respect to the width of the interval containing the jumps.  Overall, one finds that Assumption~\ref{ass.g} represents a natural generalization of the noise assumptions typically employed in noisy smooth optimization \cite{BeraByrdNoce2019,OztoByrdNoce23}, where in the smooth setting one has $\epsilon_x = 0$.
\end{example}

Let us now state our algorithm, followed by some additional justifications for its design and our previously stated assumptions. Each iteration operates in the following manner.  Consider an arbitrary iteration index $k \in \N{}$, and let $x_k \in \R{n}$ and $\epsilon_k \in \R{}_{>0}$ denote the $k$th iterate and sampling radius that have been generated by the algorithm, respectively.  The iteration begins by the algorithm sampling a set of $m \geq n+1$ points $\{x_{k,1},\dots,x_{k,m}\}$ from a uniform distribution in a ball of radius $\epsilon_k$ that is centered at~$x_k$.  (This lower bound for $m$ is consistent with GS methods for the noiseless setting; we discuss the role that it plays in our theoretical analysis in~\S\ref{sec.analysis}.)  This set of points, along with the current iterate~$x_k$, is used to construct the matrix
\bequation\label{eq.G_matrix}
  \Gtilde_k := \bbmatrix \gapprox(x_k) & \gapprox(x_{k,1}) & \cdots & \gapprox(x_{k,m}) \ebmatrix.
\eequation
Then, as in a standard gradient-sampling method, a search direction is obtained by solving the primal-dual pair of quadratic optimization problems (QPs)
\bequation\label{eq.qps}
  \left\{ \baligned \min_{(\ztilde,\dtilde) \in \R{} \times \R{n}} &\ \ztilde + \thalf \|\dtilde\|_2^2 \\ \st &\ \Gtilde_k^T\dtilde \leq \ztilde \ones \ealigned \right\}\ \ \text{and}\ \ \left\{ \baligned \max_{\ytilde \in \R{m+1}} &\ - \thalf \|\Gtilde_k\ytilde\|_2^2 \\ \st &\ \ones^T\ytilde = 1,\ \ytilde \geq 0 \ealigned \right\}.
\eequation
The final main component of an iteration of our proposed algorithm is a backtracking line search that is conducted using the approximate function values that are available.  Specifically, starting from a unit step size, the line search backtracks to attempt to obtain a step size $\alpha_k \in (0,1]$ such that, for some user-defined parameter $\eta \in (0,\thalf)$ and perturbation $\epsilonls \in \R{}_{>0}$, the following condition is satisfied:
\bequation\label{eq.armijo}
    \fapprox(x_k + \alpha_k \dtilde_k) < \fapprox(x_k) - \eta \alpha_k \|\dtilde_k\|_2^2 + \epsilonls.
\eequation
However, if \eqref{eq.armijo} fails for a trial step size $\alpha \in (0,1]$ that is short in the sense that $\alpha \|\dtilde_k\|_2 \leq \epsilon_k/3$---so that the trial point lies within a distance of $\epsilon_k/3$ of $x_k$---then the line search is abandoned, the step size is reset to zero, and the algorithm proceeds to the next iteration, where a new sample set is generated.  The motivation for this rule is that, as shown in our analysis, such a failure can only occur if the sample set is poor or the iterate is nearly stationary relative to the noise level, in which case further backtracking is not warranted. As also shown in our analysis, this procedure ensures that any step taken with a positive step size yields a decrease in the approximation function $\fapprox$ of more than $\epsilonls$, and consequently a decrease in $f$ itself.

The complete algorithm is stated as Algorithm~\ref{alg.gs}. Besides the perturbation in \eqref{eq.armijo}, there are two main differences between Algorithm~\ref{alg.gs} and a standard GS method for the noiseless setting, on which we comment in the first two bullets below. In the third bullet below, we discuss the parameter requirement in \eqref{eq.stuffff}.

\begin{algorithm}[ht]
  \caption{Noise-Tolerant Gradient Sampling Algorithm}
  \label{alg.gs}
  \begin{algorithmic}[1]
    \Require $x_1 \in \R{n}$ and $m \in \N{}$ with $m \geq n+1$
    \Require $\gamma \in (0,1)$, $\eta \in (0,\thalf)$, $\theta \in (0,1)$, $\phi \in \R{}_{>2}$, and $\nu \in \R{}_{>0}$
    \Require $\epsilon_f \in \R{}_{\geq0}$, $\fapprox : \R{n} \to \R{}$, $\epsilonls \in \R{}_{>2\epsilon_f}$, $\epsilon_x \in \R{}_{\geq0}$, $\epsilon_g \in \R{}_{>0}$, $\gapprox : \R{n} \to \R{n}$, and $\epsilon_1 \in \R{}_{>0}$ such that Assumptions~\ref{ass.f} and~\ref{ass.g} hold and
    \bequation\label{eq.stuffff}
      \epsilon_1 > \max\{\zeta,\epsilon_x\},\ \ \text{where}\ \ \zeta := \(\frac{2\epsilonls}{\eta\nu\min\{\nu,\gamma/3\}} \)^{1/2}
    \eequation
    \For{$k = 1, 2, 3, \dots$}
      \State sample $\{x_{k,1},\dots,x_{k,m}\}$ from a uniform distribution over $\ball{n}{\epsilon_k}{x_k}$
      \State set $\Gtilde_k$ by \eqref{eq.G_matrix}
      \State compute $(\ztilde_k,\dtilde_k,\ytilde_k)$ as the primal-dual solution of \eqref{eq.qps}
      \State set $\gapprox_k \gets \Gtilde_k\ytilde_k = -\dtilde_k$
      \If{$\|\gapprox_k\|_2 \leq \max\{\nu \epsilon_k, \phi\epsilon_g\}$} \label{step.nu}
        \If{$\theta \epsilon_k < \max\{\zeta,\epsilon_x\}$} \label{step.terminate}
          \State return $x_k$
        \EndIf
        \State set $\epsilon_{k+1} \gets \theta \epsilon_k$
        \State set $\alpha_k \gets 0$
      \Else
        \State set $\epsilon_{k+1} \gets \epsilon_k$
        \For{$j = 0, 1, 2, \dots$} \label{step.ls}
          \If{\eqref{eq.armijo} holds with $\alpha_k \equiv \gamma^j$}
            \State set $\alpha_k \gets \gamma^j$ and \textbf{break}
          \ElsIf{$\gamma^j \|\gapprox_k\|_2 \leq \epsilon_k/3$}\label{line.threshold}
            \State set $\alpha_k \gets 0$ and \textbf{break}
          \EndIf
        \EndFor
      \EndIf
      \State set $x_{k+1} \gets x_k + \alpha_k \dtilde_k$
    \EndFor
  \end{algorithmic}
\end{algorithm}

\bitemize[leftmargin=*]
  \item A standard GS method that employs exact objective and (generalized) gradient values takes additional steps to ensure $x_k \in \Dcal_{\nabla f}$ for all $k \in \N{}$.  We do not include such a procedure in Algorithm~\ref{alg.gs} since, in our setting, the algorithm cannot determine whether a point is included in $\Dcal_{\nabla f}$.  (Even in the noiseless setting, such a check may be impractical, which is why implementations of GS methods often skip these algorithmic steps in practice~\cite{BurkLewiOver05}.)  Our theoretical guarantees do not require such a procedure since, due to noisy function and generalized gradient evaluations, they are naturally weaker than those that can be obtained with exact values.
  \item The test in line~\ref{line.threshold} for abandoning the line search is needed due to the stochastic nature of the algorithm. For noisy smooth optimization, if an iterate is far from stationarity, then, even with a perturbed sufficient decrease condition, the line search yields a step that decreases the true objective function $f$. In Algorithm~\ref{alg.gs}, however, the line search may fail to yield sufficient decrease due merely to a bad sample set, i.e., not due to the noise. The test in line~\ref{line.threshold} identifies such an occurrence using only quantities that are available to the algorithm: if \eqref{eq.armijo} fails at a trial point that lies within a distance of $\epsilon_k/3$ of $x_k$, then either the current iterate is nearly stationary relative to the noise level or the sample set is not good in a certain sense (see~\S\ref{sec.analysis}), so the algorithm abandons the search and generates a new sample set.  The test also limits the effort that is expended in a line search, since it ensures that the number of trial step sizes in iteration $k$ is at most $\max\{0, \lceil \log(3\|\gapprox_k\|_2/\epsilon_k) / \log(1/\gamma) \rceil\} + 1$.
  \item The requirement \eqref{eq.stuffff} ensures that~$\epsilon_1$ is set appropriately for our theoretical guarantees; no such restriction is needed for a GS method in the noiseless setting. The lower bound for $\epsilon_1$ based on $\zeta$ (the bound $\epsilon_1 > \epsilon_x$ being quite natural) ensures that any step taken with a positive step size decreases $\fapprox$ by more than $\epsilonls$; see upcoming Lemma~\ref{lem.well-posed}.  Observe that $\zeta = \Ocal(\epsilonls^{1/2})$ does not depend on~$f$, and that \eqref{eq.stuffff} can be satisfied for any choices of the other parameters.  Ultimately, $\nu$ determines a trade-off between the two quantities that define the termination conditions.
\eitemize

\section{Analysis}\label{sec.analysis}

For a GS method in the noiseless setting, one can show that an algorithm similar to Algorithm~\ref{alg.gs} generates, with probability one, a sequence of iterates over which either the objective values diverge to $-\infty$ or for which any limit point is stationary for~$f$ \cite{BurkCurtLewiOverSimo20,BurkLewiOver05,CurtRobi25,Kiwi07}. Such a guarantee is not reasonable in a noisy setting, since with noise in objective and generalized-gradient evaluations, one should not expect that the algorithm would be able to generate iterates that converge all the way to stationarity for~$f$. That said, one can expect that, at any point that is far from stationarity, the algorithm will generate a search direction that leads to sufficient decrease in $\fapprox$, which in turn may correspond to sufficient decrease in $f$, at least when the norm of the step is large relative to~$\epsilon_f$. Thus, our aim in this section is to show that, from any initial point, Algorithm~\ref{alg.gs} will at least generate a sequence of iterates such that, for some $k \in \N{}$, a measure of stationarity with respect to $x_k$ is small relative to the noise in the objective function and generalized-gradient values.

Our first lemma is the following technical result. A related version of it is a hallmark of GS methods (see~\cite[Lemma 3.1]{Kiwi07}).

\blemma\label{lem.angle}
  Suppose that $\Gcal \subseteq \R{n}$ and $\tilde\Gcal \subseteq \R{n}$ are nonempty, convex, and compact sets such that with $\epsilon_g \in \R{}_{>0}$ the following conditions hold:
  \benumerate
    \item[(a)] $\dist_{\tilde\Gcal}(0) > 2\epsilon_g$;
    \item[(b)] $\dist_{\tilde\Gcal}(g) \leq \epsilon_g$ for all $g \in \Gcal$.
  \eenumerate
  Then, for any $\eta \in (0,\thalf)$, there exists $\beta \in \R{}_{>0}$ such that if $\utilde \in \tilde\Gcal$ satisfies $\|\utilde\|_2 \leq \dist_{\tilde\Gcal}(0) + \beta$, then $\utilde^T v > \eta \|\utilde\|_2^2$ for all $v \in \Gcal$.
\elemma
\bproof
  Suppose that the implication is false.  That is, suppose that for some $\eta \in (0,\thalf)$ and all $\beta \in \R{}_{>0}$, there exists $(\utilde,v) \in \tilde\Gcal \times \Gcal$ such that $\|\utilde\|_2 \leq \dist_{\tilde\Gcal}(0) + \beta$ and $\utilde^T v \leq \eta \|\utilde\|_2^2$.  Consider such $\eta \in (0,\thalf)$.  This implies that there exist infinite sequences $\{\utilde_j\}$ and $\{v_j\}$ such that for all $j \in \N{}$ one has $\utilde_j \in \tilde\Gcal$, $v_j \in \Gcal$, $\|\utilde_j\|_2 \leq \dist_{\tilde\Gcal}(0) + 1/j$, and $\utilde_j^T v_j \leq \eta \|\utilde_j\|_2^2$.  Since $\Gcal$ and $\tilde\Gcal$ are compact, the sequence $\{(\utilde_j, v_j)\}$ has a convergent subsequence.  Passing to a subsequence if necessary, one can conclude that there exists a pair of limit points $(\ubar,\vbar) \in \tilde\Gcal \times \Gcal$ such that
  \bequation\label{eq.angle}
    \ubar^T\vbar \leq \eta\|\ubar\|_2^2.
  \eequation
  On the other hand, by the definition of the sequence $\{\utilde_j\}$, it follows that $\ubar = \proj_{\tilde\Gcal}(0)$, which under condition (a) of the lemma has $\|\ubar\|_2 = \dist_{\tilde\Gcal}(0) > 2\epsilon_g$. Now let $\vhat = \proj_{\tilde\Gcal}(\vbar)$. By the Projection Theorem~\cite[Theorem~6.5]{CurtRobi25} and convexity of~$\tilde\Gcal$, one has $(0 - \ubar)^T(\vhat - \ubar) \leq 0$, i.e., $\ubar^T\vhat \geq \|\ubar\|_2^2$. Since $\vbar \in \Gcal$, condition (b) of the lemma gives $\|\vbar - \vhat\|_2 = \dist_{\tilde\Gcal}(\vbar) \leq \epsilon_g$. Thus, one finds that
  \bequationNN
    \ubar^T\vbar = \ubar^T(\vhat + \vbar - \vhat) \geq \ubar^T\vhat - \|\ubar\|_2 \|\vbar - \vhat\|_2 \geq \|\ubar\|_2^2 - \|\ubar\|_2 \epsilon_g \geq \thalf \|\ubar\|_2^2,
  \eequationNN
  where the last inequality follows since $\|\ubar\|_2>2\epsilon_g$, so $\thalf\|\ubar\|_2>\epsilon_g$, which means that $\|\ubar\|_2-\epsilon_g > \|\ubar\|_2-\thalf\|\ubar\|_2=\thalf\|\ubar\|_2$. However, this contradicts \eqref{eq.angle} since $\eta<\thalf$.
\eproof

Let us now establish some useful properties that follow by construction.  These properties hold regardless of the randomness induced by sampling.

\blemma\label{lem.well-posed}
  Suppose that Assumptions~\ref{ass.f} and \ref{ass.g} hold.  Then, for all $k \in \N{}$:
  \benumerate
    \item[(a)] $\epsilon_k \in [\max\{\zeta,\epsilon_x\},\epsilon_1]$ and $\epsilon_{k+1} \leq \epsilon_k$;
    \item[(b)] $\gapprox_k \in \gengradient_{\epsilon_x + \epsilon_k} f(x_k) + \ball{n}{\epsilon_g}{0}$;
    \item[(c)] if the line search is reached, then it tests \eqref{eq.armijo} at least once and at most a finite number of times, and it terminates with either $\alpha_k = 0$ or $\alpha_k \geq \min\{1, \frac{\gamma\epsilon_k}{3\|\gapprox_k\|_2}\}$;
    \item[(d)] if $\alpha_k > 0$, then $\fapprox(x_{k+1}) < \fapprox(x_k) - \epsilonls$ and $f(x_{k+1}) < f(x_k)$, so overall $\{\fapprox(x_k)\}$ and $\{f(x_k)\}$ are monotonically nonincreasing and each iteration with $\alpha_k > 0$ decreases $\fapprox$ by more than $\epsilonls$.
  \eenumerate
\elemma
\bproof
  Each part of the result can be proved as follows.
  \benumerate
    \item[(a)] The result follows by induction.  The algorithm requires $\epsilon_1 > \max\{\zeta,\epsilon_x\}$.  Then, $\epsilon_{k+1} = \epsilon_k$ or $\epsilon_{k+1} < \epsilon_k$, where the latter holds only if $\theta \epsilon_k \geq \max\{\zeta,\epsilon_x\}$ (otherwise, the algorithm terminated), implying $\epsilon_{k+1} = \theta\epsilon_k \geq \max\{\zeta,\epsilon_x\}$.
    \item[(b)] Each column of $\Gtilde_k$ can be expressed as $\gapprox(x)$ for some $x \in \ball{n}{\epsilon_k}{x_k}$, and by Assumption~\ref{ass.g}(a) it follows that each such column satisfies $\gapprox(x) = g(x) + e(x)$ where $g(x) \in \gengradient_{\epsilon_x} f(x) \subseteq \gengradient_{\epsilon_x + \epsilon_k} f(x_k)$ and $\|e(x)\|_2 \leq \epsilon_g$.  Hence, since $\gapprox_k$ is a convex combination of the columns of $\Gtilde_k$ and $\gengradient_{\epsilon_x + \epsilon_k} f(x_k)$ and $\ball{n}{\epsilon_g}{0}$ are convex, it follows that $\gapprox_k \in \gengradient_{\epsilon_x + \epsilon_k} f(x_k) + \ball{n}{\epsilon_g}{0}$, as claimed.
    \item[(c)] Suppose that the line search is reached.  By construction, it tests \eqref{eq.armijo} with $j = 0$, i.e., with the unit step size.  Moreover, since $\gamma \in (0,1)$ and $\epsilon_k \in \R{}_{>0}$ by part~(a), there exists a smallest $\hat\jmath \in \{0,1,2,\dots\}$ such that $\gamma^{\hat\jmath} \|\gapprox_k\|_2 \leq \epsilon_k/3$, and by construction the line search terminates with $j \leq \hat\jmath$.  Finally, if the line search terminates with $\alpha_k > 0$, then either $\alpha_k = 1$ or $\alpha_k = \gamma^j$ for some $j \geq 1$ such that the line search was not abandoned when it tested \eqref{eq.armijo} with the step size $\gamma^{j-1}$, meaning that $\gamma^{j-1} \|\gapprox_k\|_2 > \epsilon_k/3$, i.e., $\alpha_k = \gamma^j > \frac{\gamma\epsilon_k}{3\|\gapprox_k\|_2}$.
    \item[(d)] If $\alpha_k > 0$, then the line search was reached, which by construction and part~(c) means that $\|\dtilde_k\|_2 = \|\gapprox_k\|_2 > \max\{\nu \epsilon_k, \phi \epsilon_g\} \geq \nu \epsilon_k$ and $\alpha_k \geq \min\{1, \frac{\gamma\epsilon_k}{3\|\gapprox_k\|_2}\}$, so by part (a) and the definition of~$\zeta$ in the algorithm, it follows that
    \bequation\label{eq.sneaky}
      \baligned
        \eta \alpha_k \|\dtilde_k\|_2^2 &\geq \eta \min\{ \|\gapprox_k\|_2^2, \tfrac13 \gamma \epsilon_k \|\gapprox_k\|_2 \} \\
        &> \eta \nu \min\{\nu, \tfrac13\gamma\} \epsilon_k^2 \geq \eta \nu \min\{\nu, \tfrac13\gamma\} \zeta^2 = 2 \epsilonls.
      \ealigned
    \eequation
    Hence, with \eqref{eq.armijo}, one obtains $\fapprox(x_{k+1}) < \fapprox(x_k) - \eta \alpha_k \|\dtilde_k\|_2^2 + \epsilonls < \fapprox(x_k) - \epsilonls$.  Consequently, by Assumption~\ref{ass.f} and since $\epsilonls > 2\epsilon_f$, one obtains $f(x_{k+1}) \leq \fapprox(x_{k+1}) + \epsilon_f < \fapprox(x_k) - \epsilonls + \epsilon_f \leq f(x_k) - (\epsilonls - 2\epsilon_f) < f(x_k)$.
  \eenumerate
  Since each part has been proved individually, the proof is complete.
\eproof

Our next aim is to prove a lower bound for the step size obtained through the line search under certain critical conditions.  Toward this end, let us now define for any $(x,\epsilon,\beta) \in \R{n} \times \R{}_{>0} \times \R{}_{>0}$ the set
\bequationNN
  \Tcal(x, \epsilon, \beta) = \Bigg\{ S \in \prod_{i\in[m]} \ball{n}{\epsilon}{x} : \dist_{\conv(\{\gapprox(\xbar)\}_{\xbar \in S})}(0) \leq \dist_{\tilde\Gcal_{\epsilon}(x)}(0) + \beta \Bigg\}.
\eequationNN
(Here, we overload notation and use $\{\gapprox(\xbar)\}_{\xbar \in S}$ to indicate the set of values of the approximation function~$\gapprox$ evaluated at the elements in the tuple $S$.) Note that by its sampling procedure Algorithm~\ref{alg.gs} ensures for all $k \in \N{}$ that
\bequationNN
  S_k := (x_{k,1},\dots,x_{k,m}) \in \prod_{i\in[m]} \ball{n}{\epsilon_k}{x_k}.
\eequationNN
Therefore, if $S_k \in \Tcal(x_k,\epsilon_k,\beta)$, then since $\gapprox_k$ is the minimum-norm element of the convex hull of the columns of $\Gtilde_k$, it follows that
\bequation\label{eq.good_gapprox}
  \|\gapprox_k\|_2 \leq \dist_{\tilde\Gcal_{\epsilon_k}(x_k)}(0) + \beta.
\eequation

Intuitively, the set $\Tcal(x, \epsilon, \beta)$ identifies the sample sets that are good enough to guarantee that the minimum-norm element of $\tilde\Gcal_{\epsilon}(x)$ is approximated well enough by the minimum-norm element of $\conv(\{\gapprox(\xbar)\}_{\xbar \in S})$ corresponding to the sample set~$S$.  This property is important for recognizing when a point $x$ is approximately stationary.  On the other hand, if $x$ is not close to stationarity, then this property leads to a sufficiently large step size and, consequently, a sufficient decrease in the line search.

\blemma\label{lem.alpha_bound}
  Suppose Assumptions~\ref{ass.f} and \ref{ass.g} hold. Consider any $k\in\N{}$ with
  \bequation\label{eq.de}
    \dist_{\tilde\Gcal_{\epsilon_k}(x_k)}(0) > 2 \epsilon_g.
  \eequation
  Let $\eta\in(0,\thalf)$ be as in the algorithm and let $\beta_k \in \R{}_{>0}$ satisfy the conclusion of Lemma~\ref{lem.angle} with $\Gcal = \gengradient_{\epsilon_k/3} f(x_k)$ and $\tilde\Gcal = \tilde\Gcal_{\epsilon_k}(x_k)$; such a value exists since these sets are nonempty, convex, and compact (as $f$ is locally Lipschitz and Assumption~\ref{ass.g}(a) holds), \eqref{eq.de} yields Lemma~\ref{lem.angle}(a), and Assumption~\ref{ass.g}(b), which applies since $\epsilon_k \geq \epsilon_x$ by Lemma~\ref{lem.well-posed}(a), yields Lemma~\ref{lem.angle}(b). Then, if line~\ref{step.ls} is reached,
  \bequationNN
    S_k \in \Tcal(x_k, \epsilon_k, \beta_k) \implies \alpha_k \geq \min\left\{1, \frac{\gamma\epsilon_k}{3\|\gapprox_k\|_2}\right\} > 0.
  \eequationNN
\elemma
\bproof
  Under the stated conditions, suppose $S_k \in \Tcal(x_k, \epsilon_k, \beta_k)$.  By the definitions of $\tilde\Gcal_{\epsilon_k}(x_k)$ and $\Tcal(x_k, \epsilon_k, \beta_k)$ and \eqref{eq.good_gapprox}, the algorithm computes $\gapprox_k$ with $\gapprox_k \in \tilde\Gcal_{\epsilon_k}(x_k)$ and $\|\gapprox_k\|_2 \leq \dist_{\tilde\Gcal_{\epsilon_k}(x_k)}(0) + \beta_k$. Hence, by Lemma \ref{lem.angle}, it follows that
  \bequation \label{eq.lem_resul_ineq}
    \gapprox_k^Tg > \eta \|\gapprox_k\|_2^2\ \ \text{for all}\ \ g \in \gengradient_{\epsilon_k/3} f(x_k).
  \eequation
  To derive a contradiction, suppose that $\alpha_k = 0$, which by construction of the line search means it was abandoned in line~\ref{line.threshold}.  Let $\hat{\alpha}_k \in (0,1]$ be the last step size considered, for which, by construction, \eqref{eq.armijo} failed to hold and $\hat\alpha_k \|\gapprox_k\|_2 \leq \epsilon_k/3$. Failure of the line search at $\hat\alpha_k$ means by \eqref{eq.armijo} and $\dtilde_k = -\gapprox_k$ that $-\eta \hat{\alpha}_k\|\gapprox_k\|_2^2 + \epsilonls\leq\fapprox(x_k+\hat{\alpha}_k\dtilde_k)-\fapprox(x_k)$. By the noise bound in Assumption \ref{ass.f}, one thus has that
  \bequation \label{eq:ls_ineq_f_ftilde}
    -\eta \hat{\alpha}_k\|\gapprox_k\|_2^2 + \epsilonls\leq f(x_k+\hat{\alpha}_k\dtilde_k)-f(x_k)+2\epsilon_f.
  \eequation
  At the same time, Lebourg's Mean Value Theorem \cite[Theorem 5.40]{CurtRobi25} yields the existence of $\xhat_k$ on the interval $[x_k,x_k+\hat\alpha_k\dtilde_k]$ and $\ghat_k \in \gengradient f(\xhat_k)$ such that
  \bequationNN
    f(x_k+\hat\alpha_k\dtilde_k) - f(x_k) = \hat\alpha_k\ghat_k^T\dtilde_k = -\hat\alpha_k\ghat_k^T\gapprox_k,
  \eequationNN
  which with \eqref{eq:ls_ineq_f_ftilde} gives $-\eta \hat\alpha_k \|\gapprox_k\|_2^2 + \epsilonls \leq -\hat\alpha_k\ghat_k^T\gapprox_k + 2\epsilon_f$. Rearranging the terms, and using the fact that $\epsilonls>2\epsilon_f$, one obtains $\ghat_k^T\gapprox_k \leq \eta\|\gapprox_k\|_2^2 +\frac{1}{\hat\alpha_k}(2\epsilon_f-\epsilonls) \leq \eta\|\gapprox_k\|_2^2$, which with \eqref{eq.lem_resul_ineq} implies $\ghat_k \notin \gengradient_{\epsilon_k/3} f(x_k)$.  On the other hand, by the definition of the step size~$\hat\alpha_k$, one finds $\|\xhat_k - x_k\|_2 \leq \hat{\alpha}_k \|\dtilde_k\|_2 = \hat{\alpha}_k \|\gapprox_k\|_2 \leq \epsilon_k/3$, which means that $\ghat_k \in \gengradient f(\xhat_k) \subseteq \gengradient_{\epsilon_k/3} f(x_k)$, a contradiction to the conclusion of the previous sentence.  Since a contradiction has been reached, it must be true that $\alpha_k > 0$, from which the desired lower bound for $\alpha_k$ follows by Lemma~\ref{lem.well-posed}(c).
\eproof

Lemma~\ref{lem.alpha_bound} shows that, when the current iterate is not sufficiently close to stationarity, a good sample set yields a sufficiently large step size.  What remains is to ensure that good sample sets are generated sufficiently often.  For a GS method that has access to exact objective function and generalized gradient values, this follows from the sample size being large enough (specifically, $m \geq n+1$) and the objective function being continuously differentiable on an open set of full measure in $\R{n}$, which together imply that in the neighborhood of any point there exists a nonempty open set of good sample sets; see, e.g., \cite[Lemma~12.6]{CurtRobi25}.  The same argument applies in our setting if $\gapprox$ is continuous on an open set of full measure.  In Appendix~\ref{app.sampling}, we show that, in such a case, the conclusions of upcoming Lemma~\ref{lem:almost-sure-termination} and Theorem~\ref{th.main-result} hold with the assumption that we state below replaced by a lemma, as long as the sets $\tilde\Gcal_\epsilon(x)$ and $\Tcal(x,\epsilon,\beta)$ (and, correspondingly, the inclusion in Assumption~\ref{ass.g}(b)) are defined such that values of $\gapprox$ at points outside of the open set, other than at $x$ itself, are disregarded.  Here, however, so as to allow errors in the generalized gradient values that vary irregularly with $x$, we do not presume that $\gapprox$ is continuous on any open subset of $\R{n}$, and instead assume directly that good sample sets are generated with positive probability.  For consistency with the noiseless setting and with Appendix~\ref{app.sampling}, we maintain in Algorithm~\ref{alg.gs} the requirement that $m \geq n+1$, even though our assumption does not state this requirement explicitly.

We formalize our next assumption as follows.  For each $k \in \N{}$, let $\Scal_k$ represent the random tuple of sample points.  Let $(\Omega,\Fcal,\P)$ be a probability space that is generated by the sample tuples.  Now letting $\{\Fcal_k\}$ be the filtration defined by $\Fcal_1$ (defined trivially with respect to the initial conditions) and otherwise $\Fcal_k := \sigma(\Scal_1,\dots,\Scal_{k-1})$ (i.e., the $\sigma$-algebra defined by all sample tuples prior to iteration $k \geq 2$), one has, since $\fapprox$ and $\gapprox$ are measurable, that $x_k$ and $\epsilon_k$ are $\Fcal_k$-measurable while $\Scal_k$ is $\Fcal_{k+1}$-measurable.  Specifically, conditioned on $\Fcal_k$, the tuple $\Scal_k$ is uniformly distributed over $\prod_{i\in[m]} \ball{n}{\epsilon_k}{x_k}$.  We now make the following reasonable assumption.

\bassumption\label{ass.prob_good_sample}
  There exist measurable $\beta : \R{n} \times \R{}_{>0} \to \R{}_{>0}$ and $p : \R{n} \times \R{}_{>0} \to (0,1]$ such that with $\tau_\phi := (\phi+2)\epsilon_g/2$ and $\beta_\phi := (\phi-2)\epsilon_g/2$ one has:
  \benumerate
    \item[(a)] For any $(x,\epsilon) \in \R{n} \times \R{}_{>0}$ with $\epsilon \geq \epsilon_x$ and $\dist_{\tilde\Gcal_\epsilon(x)}(0) > \tau_\phi$, the value $\beta(x,\epsilon) \in (0,\beta_\phi]$ satisfies the conclusion of Lemma~\ref{lem.angle} with $\Gcal = \gengradient_{\epsilon/3} f(x)$, $\tilde\Gcal = \tilde\Gcal_{\epsilon}(x)$, and $\eta \in (0,\thalf)$ as in the algorithm.
    \item[(b)] For any $(x,\epsilon) \in \R{n} \times \R{}_{>0}$ with $\epsilon \geq \epsilon_x$ and $\dist_{\tilde\Gcal_\epsilon(x)}(0) \leq \tau_\phi$: $\beta(x,\epsilon) = \beta_\phi$.
    \item[(c)] For all $k \in \N{}$, the event $\{\Scal_k \in \Tcal(x_k,\epsilon_k,\beta(x_k,\epsilon_k))\}$ is $\Fcal_{k+1}$-measurable and
    \bequationNN
      \P(\Scal_k \in \Tcal(x_k,\epsilon_k,\beta(x_k,\epsilon_k)) | \Fcal_k) \geq p(x_k,\epsilon_k) > 0\ \ \text{almost surely}.
    \eequationNN
  \eenumerate
\eassumption

Parts~(a) and (b) of Assumption~\ref{ass.prob_good_sample} merely require a measurable selection.  Indeed, for any $(x,\epsilon)$ as in part~(a), the sets $\gengradient_{\epsilon/3} f(x)$ and $\tilde\Gcal_\epsilon(x)$ are nonempty, convex, and compact under Assumptions~\ref{ass.f} and \ref{ass.g}(a), condition~(a) of Lemma~\ref{lem.angle} holds since $\tau_\phi > 2\epsilon_g$, and condition~(b) of Lemma~\ref{lem.angle} holds by Assumption~\ref{ass.g}(b), so a value satisfying the conclusion of that lemma exists, and any smaller positive value satisfies it as well.  Hence, the substantive requirement is part~(c), which is the one that can be replaced by a lemma as shown in Appendix~\ref{app.sampling}.  Under Assumption~\ref{ass.prob_good_sample}, our next lemma shows that if the sequence of noisy objective function values is bounded below, then, with probability one, the algorithm will terminate.

\blemma\label{lem:almost-sure-termination}
  Suppose Assumptions~\ref{ass.f}, \ref{ass.g}, and \ref{ass.prob_good_sample} hold.  For any $\khat \in \N{}$, let
  \bequationNN
    \Ecal_{\khat} := \left\{\text{Algorithm~\ref{alg.gs} runs forever, $\epsilon_k = \epsilon_{\khat}$ for all $k \geq \khat$, and $\inf_{k \in \N{}} \fapprox(x_k) > -\infty$} \right\}.
  \eequationNN
  Then, $\P(\Ecal_{\khat}) = 0$. Thus, if $\displaystyle \inf_{k \in \N{}} \fapprox(x_k) > -\infty$, the algorithm terminates almost surely.
\elemma
\bproof
  Consider arbitrary $\khat \in \N{}$ generated by the algorithm.  To derive a contradiction, suppose $\P(\Ecal_{\khat}) > 0$.  On the event $\Ecal_{\khat}$, the sampling radius $\epsilon_{\khat}$ is an $\Fcal_{\khat}$-measurable random variable satisfying $\epsilon_{\khat} \geq \max\{\zeta,\epsilon_x\}$ by Lemma~\ref{lem.well-posed}(a).  Also, on the event $\Ecal_{\khat}$, for all $k \geq \khat$ the sampling radius is not reduced and the algorithm does not terminate, from which it follows that the line search is reached with
  \bequation\label{eq.hatty}
    \|\gapprox_k\|_2 > \max\{\nu \epsilon_{\khat}, \phi \epsilon_g\}\ \ \text{for all}\ \ k \geq \khat.
  \eequation
  By Lemma~\ref{lem.well-posed}(d), each iteration with $\alpha_k > 0$ decreases $\fapprox$ by more than $\epsilonls$, whereas $\alpha_k = 0$ yields $\fapprox$ unchanged.  Therefore, since $\fapprox$ is bounded below on the event $\Ecal_{\khat}$, it follows that there can only be a finite number of iterations with $k \geq \khat$ and $\alpha_k > 0$.  Hence, on the event $\Ecal_{\khat}$, there exists a random index $K \geq \khat$ such that
  \bequation\label{eq.LL}
    \alpha_k = 0\ \ \text{and}\ \ x_k = x_K\ \ \text{for all}\ \ k \geq K.
  \eequation
  Now let us define, for each $k \in \N{}$, the event $\Acal_k := \{\Scal_k \in \Tcal(x_k,\epsilon_k,\beta(x_k,\epsilon_k))\} \in \Fcal_{k+1}$. By Assumption~\ref{ass.prob_good_sample}, one has almost surely that $\P(\Acal_k | \Fcal_k) \geq p(x_k,\epsilon_k)$ with respect to the measurable function $p$.  At the same time, on the event $\Ecal_{\khat}$, one has that $(x_k,\epsilon_k) = (x_K,\epsilon_{\khat})$ for all $k \geq K$. Combined, this yields on the event $\Ecal_{\khat}$ that
  \bequationNN
    \sum_{k \in \N{}} \P(\Acal_k | \Fcal_k) \geq \sum_{k \geq K} p(x_K,\epsilon_{\khat}) = \infty.
  \eequationNN
  Thus, using the second Borel-Cantelli lemma \cite[Theorem~4.3.4]{Durr19}, it follows that almost surely on $\Ecal_{\khat}$ one has that $\Acal_k$ occurs for infinitely many~$k \geq \khat$.  Consider now any realization on $\Ecal_{\khat}$ for which $\Acal_k$ occurs infinitely often and consider any $k \geq K$ with $\Scal_k \in \Tcal(x_k,\epsilon_k,\beta(x_k,\epsilon_k))$. There are two cases to consider.
  \benumerate
    \item[(a)] If $\dist_{\tilde\Gcal_{\epsilon_{\khat}}(x_K)}(0) \leq \tau_\phi$, then by \eqref{eq.good_gapprox} and Assumption~\ref{ass.prob_good_sample}(a)--(b) (since, either way, $\beta(x_K,\epsilon_{\khat}) \leq \beta_\phi$) one has $\|\gapprox_k\|_2 \leq \tau_\phi + \beta_\phi = \phi\epsilon_g$, contradicting \eqref{eq.hatty}.
    \item[(b)] If $\dist_{\tilde\Gcal_{\epsilon_{\khat}}(x_K)}(0) > \tau_\phi > 2\epsilon_g$, then by Assumption~\ref{ass.prob_good_sample}(a) the value $\beta(x_K,\epsilon_{\khat})$ satisfies the conclusion of Lemma~\ref{lem.angle}, and since the line search is reached one has from Lemma~\ref{lem.alpha_bound} that $\alpha_k > 0$, which contradicts~\eqref{eq.LL}.
  \eenumerate
  Since a contradiction has been reached almost surely, in fact $\P(\Ecal_{\khat}) = 0$, as claimed.

  With respect to the final claim, suppose $\{\fapprox(x_k)\}$ is bounded below and the algorithm does not terminate.  Since $\epsilon_k \in [\max\{\zeta,\epsilon_x\},\epsilon_1]$ for all $k \in \N{}$ by Lemma~\ref{lem.well-posed}(a) and any reduction of the sampling radius causes it to be reduced by a factor $\theta \in (0,1)$, only a finite number of reductions can occur.  Hence, in any case, the sequence $\{\epsilon_k\}$ is eventually constant, which means that any realization of the algorithm lies in $\bigcup_{\khat \in \N{}} \Ecal_{\khat}$. However, as proved already, this is a countable union of events each of which has probability zero of occurrence, meaning that the event that $\{\fapprox(x_k)\}$ is bounded below and the algorithm does not terminate has probability zero.
\eproof

Before presenting our main theorem, we prove the following lemma that provides conditions under which the distance from the origin to a certain Goldstein generalized gradient set is bounded by the norm of the search direction computed by the algorithm, up to the noise in the generalized gradient values.

\blemma\label{eq.Gold}
  Suppose that Assumptions~\ref{ass.f} and \ref{ass.g} hold. Then, for any $k \in \N{}$,
  \bequationNN
    \dist_{\gengradient_{2\epsilon_k} f(x_k)}(0) \leq \|\gapprox_k\|_2 + \epsilon_g.
  \eequationNN
\elemma
\bproof
  By Lemma~\ref{lem.well-posed}(b), there exists $g \in \gengradient_{\epsilon_x+\epsilon_k} f(x_k)$ such that $\|g - \gapprox_k\|_2 \leq \epsilon_g$.  On the other hand, by Lemma~\ref{lem.well-posed}(a), one has $\epsilon_x \leq \max\{\zeta,\epsilon_x\} \leq \epsilon_k$, meaning that $\gengradient_{\epsilon_x+\epsilon_k} f(x_k) \subseteq \gengradient_{2\epsilon_k} f(x_k)$.  Consequently, $\dist_{\gengradient_{2\epsilon_k} f(x_k)}(0) \leq \|g\|_2 \leq \|\gapprox_k\|_2 + \epsilon_g$.
\eproof

We now present our main theorem about the behavior of Algorithm~\ref{alg.gs}.  As one can expect from state-of-the-art guarantees for GS methods in the noiseless setting, as well as from state-of-the-art guarantees for the setting of noisy smooth optimization, the theorem shows that either the iterate sequence corresponds to true objective values that decrease without bound or the algorithm will terminate and return an iterate that is approximately stationary for $f$, where the termination tolerance is only limited by the objective and generalized gradient noise bounds.

\btheorem \label{th.main-result}
  Suppose that Assumptions~\ref{ass.f}, \ref{ass.g}, and \ref{ass.prob_good_sample} hold. Then, Algorithm~\ref{alg.gs} generates sequences $\{x_k\}$ and $\{\epsilon_k\}$ such that, with probability one, either
  \begin{enumerate}[label=(\alph*)]
    \item the algorithm does not terminate and $\{f(x_k)\}\rightarrow-\infty$, or
    \item the algorithm terminates and returns $x_k \in \R{n}$ for some $k \in \N{}$ with
    \bequationNN
      \epsilon_k \in [\bar\zeta,\theta^{-1}\bar\zeta)\ \ \text{and}\ \ \|\gapprox_k\|_2 \leq \theta^{-1} \max\{\nu \bar\zeta, \phi\epsilon_g\}\ \ \text{with}\ \ \bar\zeta := \max\{\zeta,\epsilon_x\},
    \eequationNN
    from which it follows that, in the final iteration $k = k_*$, one has
    \bequationNN
      \baligned
        \epsilon_{k_*} &< \theta^{-1}\bar\zeta = \Ocal(\max\{\epsilonls^{1/2},\epsilon_x\}) \\ \text{and}\ \ 
        \dist_{\gengradient_{2\epsilon_{k_*}} f(x_{k_*})}(0) &\leq \theta^{-1} \max\{\nu \bar\zeta, \phi\epsilon_g\} + \epsilon_g = \Ocal(\epsilon_g + \max\{\epsilonls^{1/2},\epsilon_x\}),
      \ealigned
    \eequationNN
    where these bounds tend to zero as $(\epsilon_f,\epsilonls,\epsilon_x,\epsilon_g) \to (0,0,0,0)$.
  \end{enumerate}
\etheorem
\bproof
  By Lemma~\ref{lem.well-posed}(d), the sequence $\{\fapprox(x_k)\}$ is monotonically nonincreasing.  Thus, if the algorithm does not terminate, then either $\{\fapprox(x_k)\} \searrow -\infty$---meaning that (a) occurs under Assumption~\ref{ass.f}---or the sequence $\{\fapprox(x_k)\}$ is bounded below.  However, using the same logic as in the proof of the final claim of Lemma~\ref{lem:almost-sure-termination}, one can conclude that, with probability one, either the algorithm terminates or (a) occurs.

  Now suppose that the algorithm terminates in iteration $k_* \in \N{}$.  By construction, $\|\gapprox_{k_*}\|_2 \leq \max\{\nu \epsilon_{k_*}, \phi\epsilon_g\}$ and $\theta \epsilon_{k_*} < \bar\zeta$, whereas Lemma~\ref{lem.well-posed}(a) shows $\epsilon_{k_*} \geq \bar\zeta$.  Thus, $\epsilon_{k_*} \in [\bar\zeta,\theta^{-1}\bar\zeta)$, as claimed, and since $\theta \in (0,1)$ one has $\|\gapprox_{k_*}\|_2 \leq \max\{\nu \theta^{-1} \bar\zeta, \phi\epsilon_g\} \leq \theta^{-1} \max\{\nu \bar\zeta, \phi\epsilon_g\}$. The final conclusion thus follows along with Lemma~\ref{eq.Gold}.
\eproof

The dependence on $\epsilonls^{1/2}$ in Theorem~\ref{th.main-result} matches the dependence on the objective error bound that is known for noisy smooth optimization.  It arises since, by Lemma~\ref{lem.well-posed}(c), any step that is accepted with a step size less than one has $\eta\alpha_k\|\dtilde_k\|_2^2 > \tfrac13\eta\gamma\epsilon_k\|\gapprox_k\|_2$; i.e., the decrease that is guaranteed through \eqref{eq.armijo} is proportional to the product of the sampling radius and the norm of the search direction.  This is in line with the classical observation of Goldstein~\cite{Gold1977} that a step of length $\epsilon$ along the negative of the minimum-norm element of $\gengradient_\epsilon f(x)$ decreases $f$ by at least~$\epsilon$ times the norm of that element.  One consequence is that, when $f$ is smooth, our guarantee reduces to one of the type that is known for line search methods for noisy smooth optimization, as shown in the following corollary.

\bcorollary\label{cor.smooth}
  Suppose that the conditions of Theorem~\ref{th.main-result} hold, that $f$ is continuously differentiable, and that there exists $L_{\nabla f} \in \R{}_{>0}$ such that $\|\nabla f(x) - \nabla f(\xbar)\|_2 \leq L_{\nabla f} \|x - \xbar\|_2$ for all $(x,\xbar) \in \R{n} \times \R{n}$.  Then, if the algorithm terminates, the iterate $x_{k_*}$ that it returns satisfies
  \bequationNN
    \|\nabla f(x_{k_*})\|_2 \leq \theta^{-1} \max\{\nu \bar\zeta, \phi\epsilon_g\} + \epsilon_g + 2 L_{\nabla f} \theta^{-1} \bar\zeta = \Ocal(\epsilon_g + \max\{\epsilonls^{1/2},\epsilon_x\}).
  \eequationNN
\ecorollary
\bproof
  Since $f$ is continuously differentiable, one has $\gengradient f(x) = \{\nabla f(x)\}$ for all $x \in \R{n}$.  Hence, for any $\epsilon \in \R{}_{>0}$ and $x \in \R{n}$, every element of $\Ucal_{f,\epsilon}(x)$ lies in the closed and convex set $\ball{n}{L_{\nabla f}\epsilon}{\nabla f(x)}$, meaning that the same is true of every element of $\gengradient_\epsilon f(x)$, from which it follows that $\|\nabla f(x)\|_2 \leq \dist_{\gengradient_\epsilon f(x)}(0) + L_{\nabla f} \epsilon$.  The result follows from this fact with $(x,\epsilon) = (x_{k_*}, 2\epsilon_{k_*})$ and the bounds in Theorem~\ref{th.main-result}(b).
\eproof

Theorem~\ref{th.main-result} does not provide a bound on the number of iterations that are performed up to termination.  That said, if the algorithm terminates in iteration $k_* \in \N{}$, then two quantities are bounded deterministically.  By Lemma~\ref{lem.well-posed}(d), the number of iterations in which a positive step size is taken is at most $(\fapprox(x_1) - \fapprox(x_{k_*}))/\epsilonls$, and by Lemma~\ref{lem.well-posed}(a) the number of iterations in which the sampling radius is reduced is at most $\log(\epsilon_1/\bar\zeta)/\log(1/\theta)$.  The remaining iterations are those in which the line search is abandoned, the number of which is finite with probability one, but depends on the probabilities with which good sample sets are generated.  It is also worthwhile to note the role that is played by $\eta$ in the bounds of the theorem.  Since $\nu\zeta \geq (2\epsilonls/\eta)^{1/2}$ for any $\nu \in \R{}_{>0}$, the bound on $\|\gapprox_{k_*}\|_2$ is informative only if $\eta$ is not too small relative to $\epsilonls$.  Hence, while extremely small values of this parameter are common in implementations for the noiseless setting, our guarantees favor moderate values.

We close this section by noting that, in practice, if the theoretical termination tolerance of the algorithm is not known, then as in the context of noisy smooth optimization, it is possible for the algorithm to produce an iterate satisfying the conditions of Theorem~\ref{th.main-result}, then continue to produce iterates that do not satisfy such approximate stationarity conditions.  However, in such cases, since the generalized gradient values have grown in norm again, one might expect the algorithm to produce directions that reduce the true objective again, and thus to return to approximately stationary points repeatedly, although our analysis does not guarantee such behavior.

\section{Numerical Results}\label{sec.numerical}

In this section, we demonstrate the performance of Algorithm~\ref{alg.gs} on two test problems in which bounded, deterministic errors arise naturally through computational error.  The first is a fixed-order $H_\infty$  controller design problem in which each evaluation of the objective requires a search over frequency that is carried out on a lattice. The second is a PDE-constrained optimization problem for which errors arise through approximate state and adjoint solves.  In both, the bound $\epsilon_f$ on the errors in the objective values is known \textit{a~priori}, the values of $\epsilon_x$ and $\epsilon_g$ are supported by measurements, and the parameter values satisfy all of the requirements that are stated in Algorithm~\ref{alg.gs}.  Together, the experiments illustrate the behavior that is guaranteed by our analysis as well as the extent to which the guaranteed tolerances are conservative.  Two further experiments are presented in appendices: a two-variable variant of a nonsmooth Rosenbrock function with artificial noise, for which the behavior of the algorithm can be compared with quantities that are available in closed form (Appendix~\ref{app.rosenbrock}), and the training of a neural network using mini-batch estimates, in which the noise is stochastic and hence not covered by our analysis (Appendix~\ref{app.nn}).  Further details of the experiments are given in Appendix~\ref{app.details}.

\subsection{Implementation Details}\label{sec.implementation}

All of our experiments were conducted using NonOpt via its Python interface: \href{https://github.com/frankecurtis/NonOpt}{https://github.com/frankecurtis/NonOpt}. NonOpt is an open-source C++ software package for solving nonsmooth and/or nonconvex minimization problems \cite{CurtZebi25}. We employed its gradient-sampling framework through its \texttt{GradientCombination} choice for its \texttt{direction\_computation} option along with its \texttt{Backtracking} choice for its \texttt{line\_search} option.  By default, this framework employs quasi-Newton Hessian approximations, adaptive sampling, and other strategies that are not part of Algorithm~\ref{alg.gs}; hence, we configured the software to perform the iterations of Algorithm~\ref{alg.gs}.  Specifically, through the options of the software, we disabled the quasi-Newton updates (so that $\dtilde_k = -\gapprox_k$), the reuse of points from previous iterations, two heuristics by which a search direction may be accepted before any points are sampled, and the inexact solution and other modifications of the subproblems in \eqref{eq.qps} that are employed by default. We set the number of points that are sampled to $m = n+1$, and we set the line search to start from a unit step size and to employ \eqref{eq.armijo}.  Five aspects of Algorithm~\ref{alg.gs} cannot be imposed through the options of NonOpt, so we imposed them through the routines that provide the values of~$\fapprox$ and~$\gapprox$: (i) points are sampled in every iteration, whereas the software would not sample in an iteration in which a unit step along $-\gapprox(x_k)$ yields a sufficiently large reduction in the objective; (ii) the subproblem involves the generalized gradient approximations at $x_k$ and at the sample points only, whereas the software might include the one at $x_k - \gapprox(x_k)$ as well; (iii) no step is taken in an iteration in which the sampling radius is reduced; (iv) the line search is abandoned by the test in line~\ref{line.threshold}; and (v) the sampling radius is reduced, and the algorithm terminates, exactly as prescribed by the test in line~\ref{step.nu}, which differs from the test that NonOpt employs by default.

Aspect~(v) is needed since, by default, NonOpt reduces the sampling radius (or terminates) according to a test that differs from the one in line~\ref{step.nu}. Specifically, it involves $\|\cdot\|_\infty$ in place of $\|\cdot\|_2$ and an additional scaling factor, and it is tested first with $\gapprox(x_k)$ in place of $\gapprox_k$, before any points are sampled.  We disabled that test and imposed the one in line~\ref{step.nu} through the routines that provide the values of $\fapprox$; the details are provided in Appendix~\ref{app.details}.  With this arrangement, the sampling radius is reduced, and the algorithm terminates, exactly when $\|\gapprox_k\|_2 \leq \max\{\nu\epsilon_k,\phi\epsilon_g\}$ with $\gapprox_k$ computed after the $m$ points have been sampled, and never otherwise.  (Each iteration in which this test is satisfied is followed by one iteration of the software in which the sampling radius is reduced or the run terminates, and in which no step is taken. Note that we do not count this as an iteration of Algorithm~\ref{alg.gs}.)  In our experiments, the parameters are such that $\nu\max\{\zeta,\epsilon_x\} \geq \phi\epsilon_g$ for some $\phi \in \R{}_{>2}$, so the term $\phi\epsilon_g$ never decides the test, since $\epsilon_k \geq \max\{\zeta,\epsilon_x\}$ for all $k \in \N{}$.

As for the parameters of Algorithm~\ref{alg.gs}, in all of our experiments we employed $\gamma = 1/2$, $\eta = 0.45$, $\theta = 0.1$, and $\epsilon_1 = \theta^{-1}\max\{\zeta,\epsilon_x\}$.  With these choices, with $m = n+1$, and with $\epsilonls > 2\epsilon_f$, all of the requirements that are stated in Algorithm~\ref{alg.gs} are satisfied, the sampling radius is reduced exactly once---namely, to its smallest admissible value $\max\{\zeta,\epsilon_x\}$---and the algorithm terminates the first time that the test in line~\ref{step.nu} is satisfied thereafter.  The values of $\epsilon_f$, $\epsilon_x$, $\epsilon_g$, $\epsilonls$, and $\nu$ varied across our experiments, as stated in the following subsections and in Appendices~\ref{app.rosenbrock} and~\ref{app.nn}.

Our choices of $\eta$ and $\nu$ can be justified with some further explanation.  As discussed following Corollary~\ref{cor.smooth}, the guarantees in Theorem~\ref{th.main-result} favor values of $\eta$ that are not small, so we chose a value near its upper limit of $1/2$.  As for $\nu$, if $\nu \geq \gamma/3$, then \eqref{eq.stuffff} gives $\zeta = (6\epsilonls/(\eta\gamma\nu))^{1/2}$, so the smallest admissible sampling radius $\zeta$ decreases as $\nu$ increases, whereas the corresponding tolerance $\nu\zeta$ for $\|\gapprox_k\|_2$ increases with~$\nu$.  If $\nu \leq \gamma/3$, then $\nu\zeta = (2\epsilonls/\eta)^{1/2}$, which is the smallest value that this tolerance can take, whereas $\zeta$ increases as $\nu$ decreases.  Thus, a suitable value of $\nu$ depends on the scaling of the problem, and no choice of $\nu$ yields both a sampling radius and a tolerance for $\|\gapprox_k\|_2$ that are small relative to $\epsilonls^{1/2}$.

\subsection{Fixed-Order $H_{\infty}$ Controller Design}\label{sec.hinf}

A motivating setting for our algorithm is fixed-order $H_\infty$  controller design, in which each objective evaluation requires a search over frequency that is carried out only to a tolerance and that can select the wrong frequency near a point of nondifferentiability.  Consider $\dot{x} = Ax + Bu$ with $y = Cx$ under static output feedback $u = Ky$, which gives the closed-loop matrix $A_K := A + BKC$. In this context, the entries of the matrix $K \in \R{n_u \times n_y}$ are the optimization variables, so $n = n_u n_y$.  With
\bequation\label{eq.hinf.r}
  r(K) := \min_{\omega \in \R{}_{\geq0}} \sigma_{\min}(A_K - i\omega I),
\eequation
which for Hurwitz $A_K$ is its complex stability radius and equals $1/\|(zI-A_K)^{-1}\|_{H_\infty}$ \cite{Byers1988,HinrPrit1986,VanLoan1985}, let $\varsigma(K) := 1$ if $A_K$ is Hurwitz and $\varsigma(K) := -1$ otherwise, let the corresponding product be $s(K) := \varsigma(K) r(K)$, and let the objective function be
\bequation\label{eq.hinf.objective}
  f(K) := -s(K) + \thalf\rho\|K\|_F^2
\eequation
for a given $\rho \in \R{}_{>0}$, so that minimizing the function $f$ minimizes the $H_\infty$-norm of the closed-loop resolvent with a penalty on the gain.  We write the objective through the function~$s$ rather than through the $H_\infty$-norm itself because the latter is $+\infty$ for destabilizing $K$, which Assumption~\ref{ass.f} does not allow, whereas this objective function~$f$ is real-valued and locally Lipschitz everywhere. (By Weyl's inequality the function $\sigma_{\min}(\cdot - i\omega I)$ is $1$-Lipschitz in the spectral norm uniformly in $\omega$, and any segment joining a Hurwitz to a non-Hurwitz matrix meets one with an eigenvalue on the imaginary axis at which $r = 0$, so $|s(K)-s(K')| \leq \|B\|_2\|C\|_2\|K-K'\|_F$.)  The objective is nonsmooth, since two frequencies may attain the minimum in \eqref{eq.hinf.r} or the smallest singular value of $A_K - i\omega I$ at a minimizing frequency may have multiplicity greater than one, and nonconvex. Also, it is coercive and bounded below since $r(K) \leq \|A\|_2 + \|B\|_2\|C\|_2\|K\|_F$, and a minimizing frequency lies in $[0,2\|A_K\|_2]$ since $\sigma_{\min}(A_K - i\omega I) \geq \omega - \|A_K\|_2$.  Where, in the definition of $r$, the minimizing frequency $\omega^*$ is unique and the smallest singular value of $A_K - i\omega^* I$ has multiplicity one, with left and right singular vectors $u$ and $v$, one has $\nabla_K r(K) = B^T\mathrm{Re}(uv^H)C^T$.

Our instance is a lightly damped chain of three masses with springs to both walls, with forces on masses $1$ and $3$ and measurements of the position of mass $2$ and of the velocities of masses $1$ and $3$, so that $n_x = 6$, $K$ is $2 \times 3$, and $n = 6$.  The masses are $1.0$, $1.3$, and $0.80$, the spring constants are $2.0$, $1.1$, $1.4$, and $2.6$, the damping matrix is $0.015M + 0.012K_s$ with $M$ and $K_s$ the mass and stiffness matrices, and $\rho = 0.0020$; then $\|A\|_2 = 5.62$, $\|B\|_2 = 1.25$, and $\|C\|_2 = 1.00$, the open loop is Hurwitz, and $K_1 = 0$ is stabilizing with $r(K_1) = 0.0134$ and $f(K_1) = -0.0134$.

The approximations replace the minimum in \eqref{eq.hinf.r} by the minimum over a lattice $\{0,h,2h,\dots\}$ that covers $[0,2\|A_K\|_2]$.  Writing $\tilde{r}_h$ for that approximate value and $\tilde{\omega}$ for the lattice point attaining it (the smallest on ties), we set $\fapprox(K) := -\varsigma(K)\tilde{r}_h(K) + \thalf\rho\|K\|_F^2$ and $\gapprox(K) := -\varsigma(K) B^T\mathrm{Re}(uv^H)C^T + \rho K$, the latter being computed from the singular vectors at $\tilde{\omega}$.  As shown in Appendix~\ref{app.details}, for any minimizing frequency $\omega^*$ one has $\sigma_{\min}(A_K - i\omega I)^2 \leq r(K)^2 + (\omega - \omega^*)^2$ for all $\omega \in \R{}$.  Since a lattice with spacing $h$ that covers $[0,2\|A_K\|_2]$ has a point within $h/2$ of $\omega^*$, it follows that
\bequation\label{eq.hinf.epsf}
  0 \leq \tilde{r}_h(K) - r(K) \leq (r(K)^2 + h^2/4)^{1/2} - r(K).
\eequation
Given $\epsilon_f$, we first compute the lattice minimum $\tilde{r}_0(K)$ with spacing $h_0 = 10^{-2}$, for which \eqref{eq.hinf.epsf} gives $r(K) \geq r_{\rm low}(K) := (\max\{\tilde{r}_0(K)^2 - h_0^2/4,0\})^{1/2}$, and then employ
\bequation\label{eq.hinf.h}
  h(K) := \max\{2\epsilon_f,\, 2\min\{(\epsilon_f^2 + 2r_{\rm low}(K)\epsilon_f)^{1/2},\, 5 r_{\rm low}(K)\epsilon_f^{1/2}\}\}
\eequation
as the spacing. Here, the first term in the minimum is the largest spacing for which the right-hand side of \eqref{eq.hinf.epsf} is at most $\epsilon_f$ whenever $r(K) \geq r_{\rm low}(K)$, so Assumption~\ref{ass.f} holds with the given $\epsilon_f$, which is thus known \textit{a~priori}. The second term in the minimum limits the slope of $\sigma_{\min}(A_K - i\omega I)$ with respect to $\omega$ at $\tilde{\omega}$, with which the error in $\gapprox$ grows.  In our runs, $h(K)$ was between $2\epsilon_f$ and $1.5\epsilon_f^{1/2}$.  The lattices have up to $56{,}343{,}319$ points, so we evaluate each lattice minimum with a Lipschitz pruning rule (see Appendix~\ref{app.details}), which returns exactly the brute-force lattice minimum while computing at most $3{,}052$ singular value decompositions per evaluation in our runs.  The situation is akin to that illustrated on the lower-right of Figure~\ref{fig.example}: where the selected lattice point lies in the basin of $\omega^*$, the error in $\gapprox$ is of the order of $h$, whereas when two local minima over $\omega$ nearly tie the lattice may select the wrong one, which gives an error that is not small but that occurs only near a point of nondifferentiability---the behavior that $\epsilon_x$ is designed to capture.

We ran NonOpt as described in~\S\ref{sec.implementation}.  Since the generalized gradients of $f$ are small relative to the scale of the variables (their norms were between $0.049$ and $0.42$ at the points at which $\gapprox$ was evaluated), we chose $\nu = \gamma/3$, which yields the smallest possible value of $\nu\zeta$, namely, $(2\epsilonls/\eta)^{1/2}$.  All runs started from $K_1 = 0$, with $\epsilon_f \in \{10^{-5},10^{-6},10^{-7}\}$ and $\epsilonls = 2.1\epsilon_f$.  Computing $\gapprox_k$ at $K_1$ with the sampling radius $\zeta$, the largest $\epsilon_f$ for which $\|\gapprox_k\|_2 > \nu\zeta$ is approximately $0.0014$; for larger values of $\epsilon_f$ the algorithm terminates at $K_1$.  In light of the evidence presented below, we employed $\epsilon_g = 1.25\epsilon_f^{1/2}$ and $\epsilon_x = 60\epsilon_f$, which give $\max\{\zeta,\epsilon_x\} = \zeta = 18.3\epsilon_f^{1/2}$, $\epsilon_1 = 183\epsilon_f^{1/2}$, and $\nu\zeta = 3.06\epsilon_f^{1/2} \geq 2.4\epsilon_g$ at every level.  For each level we performed five runs that differ only in the seed with which the sample points are drawn.

\begin{figure}[htbp]
  \centering
  \includegraphics[width=0.8\textwidth]{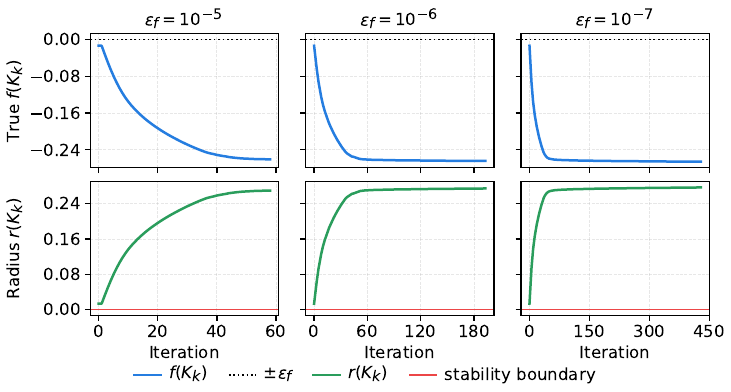}
  \caption{Fixed-order $H_\infty$  controller design: true objective values $f(K_k)$ (top) and stability radii $r(K_k)$ (bottom) over $k$ obtained by NonOpt across the three noise levels $\epsilon_f \in \{10^{-5},10^{-6},10^{-7}\}$, for the first sampling seed.  The dotted lines mark $\pm\epsilon_f$ and the solid lines mark the stability boundary $r = 0$.  The true values were computed for plotting only.}
  \label{fig.hinf.convergence}
\end{figure}

Table~\ref{tab.hinf} and Figure~\ref{fig.hinf.convergence} summarize the runs.  Every run terminated finitely through the test in line~\ref{step.terminate} of Algorithm~\ref{alg.gs}. Consistent with Lemma~\ref{lem.well-posed}(d), each of the $3{,}337$ steps with a positive step size decreased $\fapprox$ by more than $\epsilonls$ (in fact, by at least $4.3\epsilonls$) and decreased the true objective $f$.  Since $A$ is Hurwitz, $K_1 = 0$, and every accepted step decreases $f$, every iterate satisfies $f(K_k) \leq f(K_1) < 0$ and is therefore stabilizing; this was so in every run, with the stability radius improved from $r(K_1) = 0.0134$ to $0.269$--$0.270$, $0.274$, and $0.277$ for the three noise levels, respectively.

Each run terminated with $\epsilon_k = \zeta$ and $\|\gapprox_k\|_2 \leq \nu\epsilon_k$, hence within the bound of Theorem~\ref{th.main-result}(b).  Moreover, at the final iterate of every run, the minimum-norm element of the convex hull of true generalized gradients at $K_k$ and at $300$ points sampled uniformly from $\ball{n}{2\epsilon_k}{K_k}$ had norm at most $0.0061$, $0.0028$, and $6.7 \times 10^{-4}$ for the three noise levels, respectively, whereas the bound of Theorem~\ref{th.main-result}(b) on $\dist_{\gengradient_{2\epsilon_k} f(K_k)}(0)$ is $0.10$, $0.031$, and $0.010$.

\begin{table}[htbp]
  \footnotesize
  \caption{Fixed-order $H_\infty$  controller design: results obtained by NonOpt for each noise level, over five sampling seeds, reported as ranges where they differ.  Here, ``steps'' and ``abandoned'' are the numbers of iterations in which a positive step size was taken and in which the line search was abandoned by the test in line~\ref{line.threshold} of Algorithm~\ref{alg.gs}, respectively, ``evaluations'' is the total number of evaluations of $\fapprox$ and $\gapprox$, ``time'' is wall-clock time, the columns for $f$, $r$, $\|\gapprox_k\|_2$ (rounded up), and $\nu\epsilon_k$ (rounded down) pertain to the final iteration, and ``bound'' is $\theta^{-1}\nu\zeta$ in Theorem~\ref{th.main-result}(b).}
  \label{tab.hinf}
  \begin{center}
  \setlength{\tabcolsep}{5pt}
  \begin{tabular}{@{}ccccccc@{}}
    \toprule
    $\epsilon_f$ & $\zeta$ & iterations & steps & abandoned & evaluations & time (s) \\
    \midrule
    $10^{-5}$ & $0.058$ & $55$--$59$ & $53$--$57$ & $0$ & $519$--$555$ & $4$ \\
    $10^{-6}$ & $0.018$ & $194$--$202$ & $191$--$199$ & $1$--$2$ & $1{,}775$--$1{,}850$ & $27$--$30$ \\
    $10^{-7}$ & $0.0058$ & $430$--$438$ & $414$--$418$ & $12$--$18$ & $3{,}928$--$3{,}997$ & $124$--$133$ \\
    \bottomrule
  \end{tabular}

  \medskip

  \begin{tabular}{@{}cccccc@{}}
    \toprule
    $\epsilon_f$ & $f$ & $r$ & $\|\gapprox_k\|_2$ & $\nu\epsilon_k$ & bound \\
    \midrule
    $10^{-5}$ & $-0.261$--$-0.260$ & $0.269$--$0.270$ & $0.00916$--$0.00965$ & $0.00966$ & $0.096$ \\
    $10^{-6}$ & $-0.265$ & $0.274$ & $0.003009$--$0.003051$ & $0.003055$ & $0.030$ \\
    $10^{-7}$ & $-0.266$ & $0.277$ & $8.85$--$9.25 \times 10^{-4}$ & $9.66 \times 10^{-4}$ & $0.0096$ \\
    \bottomrule
  \end{tabular}
  \end{center}
\end{table}

The errors are not small relative to their bounds.  Over the iterates of each run, the largest value of $|\fapprox(K_k) - f(K_k)|$ was between $0.50\epsilon_f$ and $0.86\epsilon_f$.  As for the generalized gradients, let $g(K)$ denote the gradient of $f$ at $K$.  At the iterates of the runs and at points sampled from the balls $\ball{n}{\epsilon_k}{K_k}$ around them at which the lattice selected a frequency in the basin of the exact minimizing frequency (between $1{,}455$ and $3{,}000$ points per noise level), the error $\|\gapprox(K) - g(K)\|_2$ was at most $0.94\epsilon_f^{1/2}$.

\begin{figure}[htbp]
  \centering
  \includegraphics[width=0.8\textwidth]{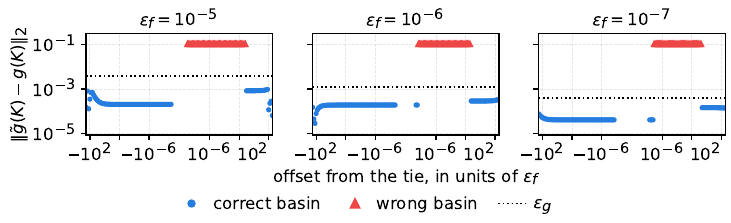}
  \caption{Fixed-order $H_\infty$  controller design: $\|\gapprox(K) - g(K)\|_2$ on both sides of the point nearest the best point found at which two frequencies tie, along the direction that separates them, against the offset from the tie in units of $\epsilon_f$, for each of the three noise levels.  Blue circles are points at which the lattice selected a frequency in the basin of the exact minimizing frequency; red triangles are those at which it did not.  The dotted line marks $\epsilon_g = 1.25\epsilon_f^{1/2}$.}
  \label{fig.hinf.kink}
\end{figure}

Figure~\ref{fig.hinf.kink} shows what happens near a point of nondifferentiability.  The best point $K^*$ that we found lies within $1.0 \times 10^{-5}$ of a point at which the two lowest local minima over $\omega$ in \eqref{eq.hinf.r} tie, and the figure shows the error in $\gapprox$ on both sides of that point along the direction that separates the two branches.  Where the lattice selects a frequency in the correct basin, the error is small, as above.  Where it selects the other basin, the error is approximately $0.11$, whatever the value of $\epsilon_f$.  However, this occurred only within $0.18\epsilon_f$ of the tie, and at every such point $\gapprox$ is within $0.20\epsilon_f^{1/2}$ of a generalized gradient of $f$ at a point that is at most $0.18\epsilon_f$ away, in agreement with Assumption~\ref{ass.g}(a).  This is the kind of error that $\epsilon_x$ captures.   By \eqref{eq.hinf.epsf} and \eqref{eq.hinf.h}, the lattice can select the other basin only where the two differ by less than $\epsilon_f$, which, to first order, is within $9.5\epsilon_f$ of the tie; our choice $\epsilon_x = 60\epsilon_f$ leaves a margin.

Finally, it is worthwhile to note that $\gapprox$ is continuous on an open set of full measure, as is required for the result in Appendix~\ref{app.sampling}; see Appendix~\ref{app.details} for details.

\subsection{PDE-Constrained Optimization with Computational Noise}
\label{sec.pde}

Our second test problem is of higher dimension, and its noise arises through the inexact solution of a PDE.  Here, Assumption~\ref{ass.f} holds with an error bound $\epsilon_f$ that is known \textit{a~priori} and that is determined by the tolerance imposed on a PDE solver.  (Algorithm~\ref{alg.gs} is applicable in settings in which the noise is uncontrollable; we consider controllable noise here since it allows for controlled experimentation.)  We also demonstrate empirically that the errors in the generalized gradients are of both of the types that are decoupled in Assumption~\ref{ass.g}, namely, typical derivative errors of size $\epsilon_g$ arising from an inexact adjoint solve, and distinct errors arising because an inexactly computed state can misidentify the correct segment in a maximum, i.e., the phenomenon that the value $\epsilon_x$ is designed to capture.

The problem that we consider, in which inexact simulation is viewed as a deterministic source of bounded error as described in~\cite{MoreWild2011}, combines ingredients that are common in the PDE-constrained optimization literature---namely, a semilinear elliptic state equation with a corresponding adjoint representation of derivatives (see, e.g., \cite{Troeltzsch2010}) and an $\ell_1$-norm penalty on the control as in \cite{CasasHerzogWachsmuth2012,Stadler2009}---with a Chebyshev ($L^\infty$) misfit, a minimum over two competing targets, and a concave peak-reward term.

Consider the two-dimensional domain $\Omega = (0,1)^2$.  Given a control $x \in \R{n}$, let the state of the system $u$ be the solution of the semilinear elliptic PDE given by
\bequation\label{eq.pde}
  -\Delta u + c u^3 = \sigma \sum_{j=1}^{n} x_j \chi_j\ \ \text{in}\ \ \Omega\ \ \text{while}\ \ u = 0\ \ \text{on}\ \ \partial\Omega,
\eequation
where $(c,\sigma) \in \R{}_{>0} \times \R{}_{>0}$ are given parameters, $\chi_j$ is the indicator function of the $j$th cell of a uniform $8 \times 8$ partition of $\Omega$, and $x$ (the optimization variable) represents a piecewise-constant array of actuators. Discretizing \eqref{eq.pde} with a standard five-point Laplacian $A$ on a $32 \times 32$ interior grid yields the discrete state equation
\bequation\label{eq.state}
  R(u,x) := Au + c u^3 - \sigma Bx = 0,
\eequation
where $A \in \R{1024 \times 1024}$, $B \in \R{1024 \times n}$, and $u \in \R{1024}$ with $n = 64$. Since $A$ is positive definite and $cu^3$ is monotone in $u$, the mapping $R(\cdot,x)$ yields unique $u = u(x)$ for each~$x$.  Moreover, since $J(u) := \partial R(u,x)/\partial u = A + 3c \diag(u^2) \succ 0$, the implicit function theorem shows that $u$ is a smooth function of~$x$ with Jacobian $\sigma J(u)^{-1} B$. Given two prescribed target profiles $u_d^{(1)}$ and $u_d^{(2)}$ and an observation index set $I$, the objective function that we consider is defined by
\bequation\label{eq.pde.objective}
  f(x) = \min_{\ell \in \{1,2\}}\ \max_{i \in I} \left| u_i(x) - u_{d,i}^{(\ell)} \right| - \mu \max_{i \in I} u_i(x) + \lambda \|x\|_1,
\eequation
where $\mu \in \R{}_{>0}$ and $\lambda \in \R{}_{>0}$ are given parameters.  That is, the objective seeks to drive the state to whichever of two profiles is easier to match in a worst-case sense while rewarding a large peak state value and penalizing non-sparse actuation. The objective is nonsmooth due to the inner $\max$, the outer $\min$, and the $\ell_1$-norm term, and it is nonconvex since a minimum of convex functions is generally nonconvex, the peak-reward term is concave in $u$, and $u$ is a nonlinear function of~$x$.  For our formulation, each target is a sum of two Gaussian bumps whose main bumps lie in opposite corners of $\Omega$, so that the two branches of the minimum have well-separated basins; see Appendix~\ref{app.details}. Considering~\eqref{eq.state} shows that $\|u(x)\|_\infty = \Ocal(\|x\|_\infty^{1/3})$, so $f$ is coercive and bounded below for any $\lambda \in \R{}_{>0}$.  We set the constants $c = 10$, $\sigma = 20$, $\lambda = 10^{-2}$, and $\mu = 1/4$ with the index set~$I$ taken to be all interior grid points in $\Omega$.

For our experiments, the approximation functions $\fapprox$ and $\gapprox$ are obtained by solving the residual equation~\eqref{eq.state} and the associated adjoint equation inexactly. Let us describe the computational procedures for computing these values, starting with $\fapprox$. The state equation is solved by a damped fixed-point (Picard) iteration, started from $u = 0$ and terminated for any given $x$ as soon as $\utilde$ is obtained satisfying $\|R(\utilde, x)\|_\infty \leq \tau$ for a prescribed residual tolerance~$\tau$.  (If the tolerance is not attained within 500 sweeps, then the solve is completed by a safeguarded Newton iteration, so that the tolerance is attained at every point at which $\fapprox$ is evaluated; this was not needed in the runs reported below.) Subtracting the state equations at $\utilde$ and $u = u(x)$ gives $R(\utilde,x) = M (\utilde - u)$ with $M = A + c \diag(\utilde^2 + \utilde u + u^2)$. Since $\utilde^2 + \utilde u + u^2 \geq 0$ pointwise, $M$ is $A$ plus a nonnegative diagonal. Both are $M$-matrices, so $0 \leq M^{-1} \leq A^{-1}$ entrywise, from which it follows that $\|\utilde - u\|_\infty \leq \|A^{-1}\|_\infty \|R(\utilde,x)\|_\infty$. The $u$-dependent part of \eqref{eq.pde.objective} is the minimum of maxima of $1$-Lipschitz functions minus $\mu$ times a maximum, so it is $(1+\mu)$-Lipschitz with respect to $\|\cdot\|_\infty$; consequently,
\bequation\label{eq.pde.epsf}
  |\fapprox(x) - f(x)| \leq (1+\mu) \|A^{-1}\|_\infty \tau\ \ \text{for all}\ \ x \in \R{n},
\eequation
where $\|A^{-1}\|_\infty = \max(A^{-1}\ones)$ is computed once. Thus, Assumption~\ref{ass.f} holds where~$\epsilon_f$ is given implicitly by the user through the choice of~$\tau$. Let us also emphasize that~$\fapprox$ is implemented to be deterministic---in the sense that repeated evaluations at the same~$x$ return the same value---although it is discontinuous since the number of solver iterations may change as $x$ varies. This is permitted by Assumption~\ref{ass.f}.

Given any $x$, a generalized gradient approximation is computed by the adjoint approach. For any $x$, letting $\ell^*$ attain the minimum in \eqref{eq.pde.objective}, $i^*$ be an index attaining the corresponding maximum, $s$ be the sign of the active residual, and $j^*$ be an index attaining the peak, one can solve $J(\utilde)p = s e_{i^*} - \mu e_{j^*}$ (approximately) to obtain $\gapprox = \sigma B^T p + \lambda \mathrm{sign}(x)$. For our experiments, this linear system is solved inexactly by the conjugate gradient method (CG), where the termination tolerance for CG yields a bound on the resulting error in $\gapprox$. Here we emphasize that the active indices $(\ell^*,i^*,j^*)$ are determined from the \textit{inexact} state $\utilde$, so near a point of nondifferentiability of $f$ they may differ from those that the true state would select. When this occurs, $\gapprox$ essentially approximates a generalized gradient at a nearby point---within a distance of $\Ocal(\epsilon_x)$---rather than at $x$ itself, and the resulting discrepancy is not reduced by tightening the solver tolerances; as we demonstrate at the end of this subsection and in Appendix~\ref{app.details}, what those tolerances control is how close to a point of nondifferentiability, and hence at what proportion of points, it occurs. Overall, the situation is akin to that illustrated on the lower-right of Figure~\ref{fig.example}.

We ran NonOpt as described in~\S\ref{sec.implementation}, with $\nu = \gamma$ rather than the value $\gamma/3$ of \S\ref{sec.hinf}, since, as discussed in \S\ref{sec.implementation}, every $\nu \leq \gamma/3$ yields the same tolerance $\nu\zeta$, and with $\nu = \gamma/3$ that tolerance was not attained within two hours at $\epsilon_f = 10^{-6}$.  All runs started from $x_1 = \ones$ at which $f(x_1) = 1.2895$.  We considered $\epsilon_f \in \{10^{-4},10^{-5},10^{-6}\}$ and in each case employed $\epsilonls = 2.1\epsilon_f$.  For each noise level we performed five runs that differ only in the seed of the random number generator. We report ranges over those runs.  For each $\epsilon_f$, the residual tolerance $\tau$ was set from \eqref{eq.pde.epsf} so that Assumption~\ref{ass.f} holds with that value of $\epsilon_f$ (which gives $\tau \approx 11\epsilon_f$), and the CG termination tolerance for the adjoint solve was set so that the corresponding bound on the error in $\gapprox$ was equal to $\epsilon_f$.  The inexact state also perturbs the adjoint operator, since the system that is solved involves $J(\utilde)$ rather than $J(u(x))$. Since $\|J(u)^{-1}\|_2 \leq \|A^{-1}\|_2$ for all $u$ and by the bound on $\|\utilde - u\|_\infty$ above, this contributes an error in $\gapprox$ of at most $3c\sigma(1+\mu)\|B\|_2\|A^{-1}\|_2^2\|A^{-1}\|_\infty(\|\utilde\|_\infty + \|u(x)\|_\infty)\tau$, which is proportional to $\epsilon_f$.  The sum of these two bounds, which is valid wherever the active indices are identified correctly, was at most $17.8\epsilon_f$ at the points at which~$\gapprox$ was evaluated in our runs.  Accordingly, we employed $\epsilon_g = 20\epsilon_f$ and, in light of the evidence presented at the end of this subsection, $\epsilon_x = \epsilon_f$.  With these values, $\max\{\zeta,\epsilon_x\} = \zeta = 10.6\epsilon_f^{1/2}$, $\epsilon_1 = 106\epsilon_f^{1/2}$, and $\nu\zeta = 5.29\epsilon_f^{1/2} \geq 3\epsilon_g$ for all of the noise levels that we consider.  (The values $\epsilon_x$ and $\epsilon_g$ affect Algorithm~\ref{alg.gs} only through $\max\{\zeta,\epsilon_x\}$ and $\max\{\nu\epsilon_k,\phi\epsilon_g\}$, so the runs would be identical for any $\epsilon_x \leq \zeta$ and any $\epsilon_g \leq \nu\zeta/\phi$.)  For presenting our results, we employed as a reference the best objective value that we found, namely, $f_{\text{ref}} = 0.08334$.  It is the smallest true objective value over the iterates of all of the runs that we performed for this problem.  (These include runs with other parameter values that were stopped at time limits; see Appendix~\ref{app.details}.)

\begin{figure}[ht]
  \centering
  \includegraphics[width=0.8\textwidth]{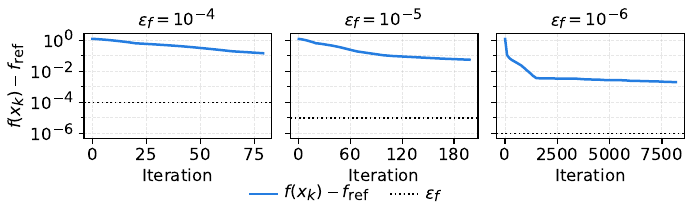}
  \caption{PDE-constrained problem: true objective value at the iterates, shown as the gap $f(x_k) - f_{\text{ref}}$, over $k$ obtained by NonOpt across the three noise levels $\epsilon_f \in \{10^{-4},10^{-5},10^{-6}\}$ with $\epsilon_{\text{ls}} = 2.1\epsilon_f$.  The dotted line in each panel marks $\epsilon_f$.  The true objective values were computed for plotting only and were never available to the algorithm.}
  \label{fig.pde.convergence}
\end{figure}

Table~\ref{tab.pde} and Figure~\ref{fig.pde.convergence} summarize the runs.  Every run terminated finitely through the test in line~\ref{step.terminate} of Algorithm~\ref{alg.gs}.  Consistent with Lemma~\ref{lem.well-posed}(d), each of the $9{,}278$ steps with a positive step size that were taken over all of the runs decreased $\fapprox$ by more than $\epsilonls$ (in fact, by at least $1.7\epsilonls$) and decreased the true objective $f$, so the true objective values in Figure~\ref{fig.pde.convergence} decrease monotonically.  The gap $f(x_k) - f_{\text{ref}}$ was reduced from $1.206$ at the initial point to $0.14$, $0.053$--$0.056$, and $0.0020$ for the three noise levels, respectively.  (Theorem~\ref{th.main-result} pertains to a measure of stationarity, not to the objective gap, so this observation complements, rather than follows from, our theoretical guarantees.)  As for stationarity, each run terminated with $\epsilon_k = \zeta$ and with a value of $\|\gapprox_k\|_2$ that is within the bound of Theorem~\ref{th.main-result}(b).  The last column of Table~\ref{tab.pde} reports the norm of the minimum-norm element of the convex hull of the true generalized gradients at the final iterate $x_k$ and at $300$ points sampled uniformly from $\ball{n}{2\epsilon_k}{x_k}$, which is an upper bound for $\dist_{\gengradient_{2\epsilon_k} f(x_k)}(0)$. It is below the corresponding bound $\theta^{-1}\nu\zeta + \epsilon_g$ of the theorem in each case, and it is smaller than the same quantity at the initial point $x_1$ by factors of at least $3.3$, $9.8$, and $36$, respectively.  For $\epsilon_f = 10^{-6}$, the line search was abandoned by the test in line~\ref{line.threshold} in between $67\%$ and $81\%$ of the iterations, none of which were among the first $555$; i.e., in the latter part of these runs, most of the sample sets yielded neither a step that satisfies \eqref{eq.armijo} nor a value of $\|\gapprox_k\|$ that satisfies the termination test, and the algorithm proceeded by generating new sample sets, as intended.  The measured errors confirm that \eqref{eq.pde.epsf} is valid, and also show its conservativeness. Over the iterates of the runs, $\max_k |\fapprox(x_k) - f(x_k)|$ was smaller than $\epsilon_f$ by a factor of between $58$ and $69$.  Since $\zeta$ is proportional to $\epsilon_f^{1/2}$, a bound that was tighter by the factor $58$ would have given sampling radii---and stationarity tolerances---smaller by a factor of about $7.6$.  For comparison, we also ran the gradient sampling method of NonOpt without the relaxation---i.e., with $\epsilonls = 0$, with no abandonment of the line search, and with its default test for reducing the sampling radius, which was allowed to decrease to $10^{-4}\zeta$---with the same $\fapprox$ and $\gapprox$ for as many evaluations as the run of Algorithm~\ref{alg.gs} with the first sampling seed.  In each case it had not terminated when these evaluations were spent and $f - f_{\text{ref}}$ was still decreasing, having reached $0.24$, $0.10$, and $0.0066$ for the three noise levels, respectively, compared with $0.14$, $0.055$, and $0.0020$ for that run of Algorithm~\ref{alg.gs}.

\begin{table}[htbp]
  \footnotesize
  \caption{PDE-constrained problem: results obtained by NonOpt for each noise level, over five sampling seeds, reported as ranges where they differ.  The columns are as in Table~\ref{tab.hinf}, the column for $f - f_{\text{ref}}$ pertains to the final iteration, and ``dist.'' is an upper bound for $\dist_{\gengradient_{2\epsilon_k} f(x_k)}(0)$ at the final iterate that was computed from true generalized gradients.}
  \label{tab.pde}
  \begin{center}
  \setlength{\tabcolsep}{5pt}
  \begin{tabular}{@{}ccccccc@{}}
    \toprule
    $\epsilon_f$ & $\zeta$ & iterations & steps & abandoned & evaluations & time (s) \\
    \midrule
    $10^{-4}$ & $0.11$ & $80$ & $78$ & $0$ & $5{,}558$ & $21$--$25$ \\
    $10^{-5}$ & $0.033$ & $195$--$210$ & $193$--$208$ & $0$ & $13{,}270$--$14{,}276$ & $63$--$68$ \\
    $10^{-6}$ & $0.011$ & $4{,}777$--$8{,}168$ & $1{,}572$--$1{,}590$ & $3{,}199$--$6{,}594$ & $325{,}152$--$555{,}157$ & $2{,}137$--$3{,}745$ \\
    \bottomrule
  \end{tabular}

  \medskip

  \begin{tabular}{@{}cccccc@{}}
    \toprule
    $\epsilon_f$ & $f - f_{\text{ref}}$ & $\|\gapprox_k\|_2$ & $\nu\epsilon_k$ & bound & dist. \\
    \midrule
    $10^{-4}$ & $0.14$ & $0.0509$--$0.0513$ & $0.0529$ & $0.52$ & $0.047$ \\
    $10^{-5}$ & $0.053$--$0.056$ & $0.01655$--$0.01668$ & $0.01673$ & $0.16$ & $0.015$--$0.016$ \\
    $10^{-6}$ & $0.0020$ & $0.00503$--$0.00527$ & $0.00529$ & $0.052$ & $0.0042$--$0.0044$ \\
    \bottomrule
  \end{tabular}
  \end{center}
\end{table}

Finally, let us illustrate the manner in which Assumption~\ref{ass.g} distinguishes between errors captured by $\epsilon_x$ and by $\epsilon_g$.  Let $x_{\text{end}}$ denote the final iterate of the run with $\epsilon_f = 10^{-6}$ for the first sampling seed.  Within a ball of radius $10^{-3}$ centered at $x_{\text{end}}$, which is a region in which $f$ is nonsmooth, we located $24$ pairs of points at which the exact state selects different active indices, bisected each segment to the point at which the selection changes, and sampled finely around it.  At every point at which the inexact state selected different indices, we checked Assumption~\ref{ass.g}(a) as it is stated: there is a point $\bar{x}$ with $\|\bar{x}-x\|_2 \leq 0.33\epsilon_f$ at which the indices selected at $x$ truly are the active ones, and an element of $\gengradient f(\bar{x})$ that differs from $\gapprox(x)$ by at most $0.052\epsilon_f$, $0.040\epsilon_f$, and $0.018\epsilon_f$ for $\epsilon_f = 10^{-4}$, $\epsilon_f = 10^{-5}$, and $\epsilon_f = 10^{-6}$, respectively, which is within $\epsilon_g$ in each case.  (At points sampled uniformly from this ball, the error $\|\gapprox(x) - g(x)\|_2$, where $g(x)$ denotes the generalized gradient obtained by the same adjoint procedure with the state and adjoint equations solved to high accuracy, was between $0.016$ and $0.043$ at every noise level where the active indices were misidentified and at most $0.078\epsilon_f$ elsewhere; see Appendix~\ref{app.details}.)  A least-squares fit of $\log_{10}$ of the largest such distance against $\log_{10}\epsilon_f$ has slope $1.0$, so the distance is proportional to~$\epsilon_f$.  This is why the decoupling in Assumption~\ref{ass.g} is critical here. Indeed, the discrepancy at such a point is not reduced by tightening the solver tolerances, whereas the distance to a point at which $\gapprox(x)$ is an accurate generalized gradient shrinks in proportion to~$\epsilon_f$ and is therefore captured by~$\epsilon_x$.  Accordingly we set $\epsilon_x = \epsilon_f$, which leaves a margin over the measurement, and $\epsilon_x \leq \zeta$ at every level.  (In the runs for Table~\ref{tab.pde}, the active indices were identified correctly at every iterate.)

Let us note that, as in \S\ref{sec.hinf}, $\gapprox$ is continuous on an open set of full measure, barring a degenerate situation (see Appendix~\ref{app.details}), as is required for the result in Appendix~\ref{app.sampling}.

\section{Conclusion}\label{sec.conclusion}

We have proposed, analyzed, and tested an algorithm for solving nonconvex and nonsmooth optimization problems when objective and generalized-gradient evaluations may be corrupted by bounded, uncontrollable errors.  To the best of our knowledge, ours is the first gradient-sampling-based method for solving nonconvex problems in such a setting. We have proved reasonable convergence guarantees for our proposed algorithm, and demonstrated through multiple experiments that the algorithm behaves as predicted by our analysis in the presence of noise.

Our numerical experiments employed a software implementation that was configured to perform the iterations of Algorithm~\ref{alg.gs}.  In the experiments with deterministic and bounded noise, the results are consistent with our analysis. The algorithm terminated finitely, every step that was taken reduced the true objective, and the points that were returned are approximately stationary to within the guaranteed tolerances.  The results also show that these tolerances are conservative.  One might expect strategies that are not covered by our analysis---e.g., quasi-Newton Hessian approximations and adaptive sampling as in \cite{CurtQue13,CurtQue15}---to yield better practical performance, as they do in the noiseless setting. In any case, in our experience, the most important parameter for such an algorithm in a noisy setting is $\epsilonls$, which depends on the noise level in the objective values, for which estimation techniques are well known~\cite{MoreWild2011}.

%% file: appendix.tex
\section{Replacing Assumption~\ref{ass.prob_good_sample} by a Lemma}\label{app.sampling}

In this appendix, we show that Assumption~\ref{ass.prob_good_sample} can be replaced by a lemma if the generalized gradient approximation function $\gapprox$ is continuous over an open set of full measure in $\R{n}$.  This is a counterpart of an assumption that underlies analyses of GS methods in the noiseless setting, namely, that $f$ is continuously differentiable over such a set~\cite{BurkLewiOver05,CurtRobi25,Kiwi07}.  It allows $\gapprox$ to have an arbitrary number of discontinuities as long as they are confined to a set of measure zero, which is typically the situation when discontinuities arise only where the identity of an active piece of $f$, or the number of iterations performed by an iterative solver, changes.  In such a setting, the property that is required by Assumption~\ref{ass.prob_good_sample}(c) can be proved as long as the sets that are involved are defined in a manner that accounts for the fact that, with probability one, the algorithm never evaluates $\gapprox$ at a sample point outside of the open set.

Given open $\Dcal_{\gapprox} \subseteq \R{n}$, let us define for any $(x,\epsilon,\beta) \in \R{n} \times \R{}_{>0} \times \R{}_{>0}$ the sets
\bequationNN
  \baligned
  \hat\Gcal_\epsilon(x) &:= \cl(\conv(\{\gapprox(x)\} \cup \{\gapprox(\xbar) : \xbar \in \ball{n}{\epsilon}{x} \cap \Dcal_{\gapprox}\})) \\
  \text{and}\ \ 
  \hat\Tcal(x,\epsilon,\beta) &:= \Bigg\{S \in \prod_{i\in[m]} (\ball{n}{\epsilon}{x} \cap \Dcal_{\gapprox}) : \\
  &\qquad\quad \dist_{\conv(\{\gapprox(x)\} \cup \{\gapprox(\xbar)\}_{\xbar \in S})}(0) \leq \dist_{\hat\Gcal_{\epsilon}(x)}(0) + \beta \Bigg\}.
  \ealigned
\eequationNN
These are counterparts of $\tilde\Gcal_\epsilon(x)$ from \S\ref{sec.algorithm} and $\Tcal(x,\epsilon,\beta)$ from \S\ref{sec.analysis}, respectively, in which values of $\gapprox$ at points outside of $\Dcal_{\gapprox}$, other than at $x$ itself, are disregarded.  In particular, one has $\hat\Gcal_\epsilon(x) \subseteq \tilde\Gcal_\epsilon(x)$, where equality holds unless values of $\gapprox$ on a set of measure zero enlarge $\tilde\Gcal_\epsilon(x)$.  We make the following assumption.

\bassumption\label{ass.g_continuous}
  There exists open $\Dcal_{\gapprox} \subseteq \R{n}$ of full measure in $\R{n}$ such that $\gapprox$ is continuous over $\Dcal_{\gapprox}$ and such that Assumption~\ref{ass.g}(b) holds with $\hat\Gcal_\epsilon(x)$ in place of $\tilde\Gcal_\epsilon(x)$, i.e., $\gengradient_{\epsilon/3} f(x) \subseteq \hat\Gcal_\epsilon(x) + \ball{n}{\epsilon_g}{0}$ for all $x \in \R{n}$ and $\epsilon \in \R{}_{>0}$ with $\epsilon \geq \epsilon_x$.
\eassumption

The inclusion in Assumption~\ref{ass.g_continuous} is stronger than the one in Assumption~\ref{ass.g}(b), but the difference is slight.  For instance, the inclusion holds in each of the examples in \S\ref{sec.algorithm} when $\Dcal_{\gapprox}$ is the complement of the set of points at which $\gapprox$ jumps and $\gapprox$ is continuous elsewhere, since the arguments in those examples only involve the value of $\gapprox$ at $x$ and values of $\gapprox$ at points on either side of a jump.  Assumption~\ref{ass.g_continuous} also holds in the noiseless setting: If $\Dcal_{\gapprox}$ is open with full measure in $\R{n}$, and over it $f$ is continuously differentiable, and if $\gapprox(x) = \nabla f(x)$ for all $x \in \Dcal_{\gapprox}$ and $\gapprox(x) \in \gengradient f(x)$ otherwise, then Assumptions~\ref{ass.g}(a) and \ref{ass.g_continuous} hold with $\epsilon_x = \epsilon_g = 0$; see \cite[Theorem~2.5.1]{Clar83}.

Let us now state a counterpart of Lemma~\ref{lem.alpha_bound}.

\blemma\label{lem.alpha_bound_hat}
  Suppose Assumptions~\ref{ass.f}, \ref{ass.g}, and \ref{ass.g_continuous} hold. Consider any $k\in\N{}$ with $\dist_{\hat\Gcal_{\epsilon_k}(x_k)}(0) > 2 \epsilon_g$.  Let $\eta\in(0,\thalf)$ be as in the algorithm and let $\beta_k \in \R{}_{>0}$ satisfy the conclusion of Lemma~\ref{lem.angle} with $\Gcal = \gengradient_{\epsilon_k/3} f(x_k)$ and $\tilde\Gcal = \hat\Gcal_{\epsilon_k}(x_k)$; such a value exists since these sets are nonempty, convex, and compact, the stated lower bound for $\dist_{\hat\Gcal_{\epsilon_k}(x_k)}(0)$ yields Lemma~\ref{lem.angle}(a), and Assumption~\ref{ass.g_continuous}, which applies since $\epsilon_k \geq \epsilon_x$ by Lemma~\ref{lem.well-posed}(a), yields Lemma~\ref{lem.angle}(b). Then, if line~\ref{step.ls} is reached,
  \bequationNN
    S_k \in \hat\Tcal(x_k, \epsilon_k, \beta_k) \implies \alpha_k \geq \min\left\{1, \frac{\gamma\epsilon_k}{3\|\gapprox_k\|_2}\right\} > 0.
  \eequationNN
\elemma
\bproof
  Under the stated conditions, suppose $S_k \in \hat\Tcal(x_k, \epsilon_k, \beta_k)$.  Then, each column of $\Gtilde_k$ is either $\gapprox(x_k)$ or the value of $\gapprox$ at an element of $\ball{n}{\epsilon_k}{x_k} \cap \Dcal_{\gapprox}$, so since $\gapprox_k$ is the minimum-norm element of the convex hull of the columns of $\Gtilde_k$, it follows that $\gapprox_k \in \hat\Gcal_{\epsilon_k}(x_k)$ and $\|\gapprox_k\|_2 \leq \dist_{\hat\Gcal_{\epsilon_k}(x_k)}(0) + \beta_k$.  Hence, by Lemma~\ref{lem.angle}, one finds that \eqref{eq.lem_resul_ineq} holds.  The remainder of the proof is identical to that of Lemma~\ref{lem.alpha_bound}, which, after \eqref{eq.lem_resul_ineq} has been established, does not involve the set $\tilde\Gcal_{\epsilon_k}(x_k)$.
\eproof

Our next lemma shows that, as in the noiseless setting (see~\cite[Lemma~12.6]{CurtRobi25}), there is a nonempty open set of good sample tuples.

\blemma\label{lem.good_samples}
  Suppose that Assumptions~\ref{ass.f}, \ref{ass.g}, and \ref{ass.g_continuous} hold and that $m \geq n+1$.  Then, for any $(x,\epsilon,\beta) \in \R{n} \times \R{}_{>0} \times \R{}_{>0}$, the set $\hat\Tcal(x,\epsilon,\beta)$ contains a nonempty open subset of $\prod_{i\in[m]} \ball{n}{\epsilon}{x}$.
\elemma
\bproof
  Consider arbitrary $(x,\epsilon,\beta) \in \R{n} \times \R{}_{>0} \times \R{}_{>0}$ and let
  \bequationNN
    \Ccal := \{\gapprox(x)\} \cup \{\gapprox(\xbar) : \xbar \in \ball{n}{\epsilon}{x} \cap \Dcal_{\gapprox}\},\ \ \text{so that}\ \ \hat\Gcal_\epsilon(x) = \cl(\conv(\Ccal)).
  \eequationNN
  By the definition of $\hat\Gcal_\epsilon(x)$, there exists $u \in \conv(\Ccal)$ with $\|u\|_2 \leq \dist_{\hat\Gcal_\epsilon(x)}(0) + \beta/2$.  Hence, by Carath\'eodory's theorem, there exist an integer $r \in \{0,1,\dots,n+1\}$, points $\{\xbar_1,\dots,\xbar_r\} \subset \ball{n}{\epsilon}{x} \cap \Dcal_{\gapprox}$, and nonnegative scalars $\{\lambda_0,\lambda_1,\dots,\lambda_r\}$ that sum to one such that $u = \lambda_0 \gapprox(x) + \sum_{i=1}^r \lambda_i \gapprox(\xbar_i)$.  Since $\Dcal_{\gapprox}$ is open and $\gapprox$ is continuous over $\Dcal_{\gapprox}$, for each $i \in [r]$ there exists $\delta_i \in \R{}_{>0}$ such that $\ball{n}{\delta_i}{\xbar_i} \subseteq \Dcal_{\gapprox}$ and $\|\gapprox(\xhat) - \gapprox(\xbar_i)\|_2 \leq \beta/2$ for all $\xhat \in \ball{n}{\delta_i}{\xbar_i}$.  For each $i \in [r]$, let $\Ucal_i$ be the intersection of the interiors of $\ball{n}{\delta_i}{\xbar_i}$ and $\ball{n}{\epsilon}{x}$, which is open and, since $\xbar_i \in \ball{n}{\epsilon}{x}$, nonempty.  For each $i \in [m] \setminus [r]$, where $m \geq n+1 \geq r$, let $\Ucal_i$ be the intersection of $\Dcal_{\gapprox}$ and the interior of $\ball{n}{\epsilon}{x}$, which is open and, since $\Dcal_{\gapprox}$ has full measure in $\R{n}$, nonempty.  Then, $\prod_{i\in[m]} \Ucal_i$ is a nonempty open subset of $\prod_{i\in[m]} (\ball{n}{\epsilon}{x} \cap \Dcal_{\gapprox})$.  Now consider arbitrary $S = (\xhat_1,\dots,\xhat_m) \in \prod_{i\in[m]} \Ucal_i$.  Since $\uhat := \lambda_0 \gapprox(x) + \sum_{i=1}^r \lambda_i \gapprox(\xhat_i)$ is an element of $\conv(\{\gapprox(x)\} \cup \{\gapprox(\xbar)\}_{\xbar \in S})$ with $\|\uhat - u\|_2 \leq \sum_{i=1}^r \lambda_i \|\gapprox(\xhat_i) - \gapprox(\xbar_i)\|_2 \leq \beta/2$, one finds
  \bequationNN
    \dist_{\conv(\{\gapprox(x)\} \cup \{\gapprox(\xbar)\}_{\xbar \in S})}(0) \leq \|\uhat\|_2 \leq \|u\|_2 + \beta/2 \leq \dist_{\hat\Gcal_\epsilon(x)}(0) + \beta,
  \eequationNN
  meaning that $S \in \hat\Tcal(x,\epsilon,\beta)$, which completes the proof.
\eproof

We now show that the conclusions of our main results hold with Assumption~\ref{ass.g_continuous} in place of Assumption~\ref{ass.prob_good_sample}.

\btheorem\label{th.appendix}
  Suppose that Assumptions~\ref{ass.f}, \ref{ass.g}, and \ref{ass.g_continuous} hold.  Then, the conclusions of Lemma~\ref{lem:almost-sure-termination} and Theorem~\ref{th.main-result} hold.
\etheorem
\bproof
  Let the probability space $(\Omega,\Fcal,\P)$, the filtration $\{\Fcal_k\}$, and the random tuples $\{\Scal_k\}$ be defined as in \S\ref{sec.analysis}, so that, conditioned on $\Fcal_k$, the tuple $\Scal_k$ is uniformly distributed over $\prod_{i\in[m]} \ball{n}{\epsilon_k}{x_k}$.  For any $\beta \in \R{}_{>0}$, let $q_\beta : \R{n} \times \R{}_{>0} \to [0,1]$ be defined such that, for any $(x,\epsilon) \in \R{n} \times \R{}_{>0}$, the value $q_\beta(x,\epsilon)$ is the probability that a tuple that is uniformly distributed over $\prod_{i\in[m]} \ball{n}{\epsilon}{x}$ is an element of $\hat\Tcal(x,\epsilon,\beta)$.  Since any nonempty open set has positive measure, Lemma~\ref{lem.good_samples} shows that $q_\beta(x,\epsilon) > 0$ for all $(x,\epsilon) \in \R{n} \times \R{}_{>0}$.  Moreover, for all $k \in \N{}$ and $\beta \in \R{}_{>0}$, the event $\{\Scal_k \in \hat\Tcal(x_k,\epsilon_k,\beta)\}$ is $\Fcal_{k+1}$-measurable and
  \bequation\label{eq.prob_good_sample}
    \P(\Scal_k \in \hat\Tcal(x_k,\epsilon_k,\beta) | \Fcal_k) = q_\beta(x_k,\epsilon_k) > 0\ \ \text{almost surely}.
  \eequation
  (The measurability properties required here hold since $\gapprox$ is measurable, $\Dcal_{\gapprox}$ is open, and, by Carath\'eodory's theorem and continuity of $\gapprox$ over $\Dcal_{\gapprox}$, the value $\dist_{\hat\Gcal_\epsilon(x)}(0)$ is the infimum over a countable family of measurable functions of $(x,\epsilon)$, namely, the norms of the convex combinations with rational weights of $\gapprox(x)$ and values of $\gapprox$ at points of the form $x + \epsilon z \in \Dcal_{\gapprox}$ with $z$ in a countable dense subset of the unit ball.)

  Now consider arbitrary $\khat \in \N{}$ along with the event $\Ecal_{\khat}$ defined in Lemma~\ref{lem:almost-sure-termination} and, to derive a contradiction, suppose that $\P(\Ecal_{\khat}) > 0$.  The arguments in the proof of Lemma~\ref{lem:almost-sure-termination} up to and including \eqref{eq.LL} do not involve Assumption~\ref{ass.prob_good_sample}.  Hence, on the event $\Ecal_{\khat}$, one has that \eqref{eq.hatty} holds and there exists a random index $K \geq \khat$ such that \eqref{eq.LL} holds.  Now let $\tau_\phi := (\phi+2)\epsilon_g/2$ and $\beta_\phi := (\phi-2)\epsilon_g/2$, where $\tau_\phi > 2\epsilon_g$ and $\beta_\phi > 0$ since $\phi > 2$ and $\epsilon_g > 0$, and let us define, for each $(k,i) \in \N{} \times \N{}$, the event $\Acal_{k,i} := \{\Scal_k \in \hat\Tcal(x_k,\epsilon_k,\beta_\phi/i)\} \in \Fcal_{k+1}$. By \eqref{eq.prob_good_sample}, one has almost surely that $\P(\Acal_{k,i} | \Fcal_k) = q_{\beta_\phi/i}(x_k,\epsilon_k)$.  At the same time, on the event $\Ecal_{\khat}$, one has that $(x_k,\epsilon_k) = (x_K,\epsilon_{\khat})$ for all $k \geq K$. Combined, this yields on $\Ecal_{\khat}$ that, for all $i \in \N{}$,
  \bequationNN
    \sum_{k \in \N{}} \P(\Acal_{k,i} | \Fcal_k) \geq \sum_{k \geq K} q_{\beta_\phi/i}(x_K,\epsilon_{\khat}) = \infty.
  \eequationNN
  Thus, using the second Borel-Cantelli lemma \cite[Theorem~4.3.4]{Durr19}, it follows for each $i \in \N{}$ that almost surely on $\Ecal_{\khat}$ one has that $\Acal_{k,i}$ occurs for infinitely many~$k \geq \khat$.  Hence, since $\N{}$ is countable, almost surely on $\Ecal_{\khat}$ one has for all $i \in \N{}$ that $\Acal_{k,i}$ occurs for infinitely many~$k \geq \khat$.  Consider now any realization on $\Ecal_{\khat}$ for which this is the case. There are two cases to consider.
  \benumerate
    \item[(a)] If $\dist_{\hat\Gcal_{\epsilon_{\khat}}(x_K)}(0) \leq \tau_\phi$, then consider any $k \geq K$ with $\Scal_k \in \hat\Tcal(x_k,\epsilon_k,\beta_\phi)$.  Since $\gapprox_k$ is the minimum-norm element of the convex hull of the columns of $\Gtilde_k$, one has $\|\gapprox_k\|_2 \leq \dist_{\hat\Gcal_{\epsilon_{\khat}}(x_K)}(0) + \beta_\phi \leq \tau_\phi + \beta_\phi = \phi\epsilon_g$, contradicting \eqref{eq.hatty}.
    \item[(b)] If $\dist_{\hat\Gcal_{\epsilon_{\khat}}(x_K)}(0) > \tau_\phi > 2\epsilon_g$, then, as argued in the statement of Lemma~\ref{lem.alpha_bound_hat}, there exists $\beta \in \R{}_{>0}$ satisfying the conclusion of Lemma~\ref{lem.angle} with $\Gcal = \gengradient_{\epsilon_{\khat}/3} f(x_K)$ and $\tilde\Gcal = \hat\Gcal_{\epsilon_{\khat}}(x_K)$, and any smaller positive value satisfies it as well.  Hence, considering any $i \in \N{}$ with $\beta_\phi/i \leq \beta$ and any $k \geq K$ with $\Scal_k \in \hat\Tcal(x_k,\epsilon_k,\beta_\phi/i)$, since the line search is reached one has from Lemma~\ref{lem.alpha_bound_hat} that $\alpha_k > 0$, which contradicts~\eqref{eq.LL}.
  \eenumerate
  Since a contradiction has been reached almost surely, it follows that $\P(\Ecal_{\khat}) = 0$.  The remainder of the proof of Lemma~\ref{lem:almost-sure-termination}, as well as the proof of Theorem~\ref{th.main-result}, do not involve Assumption~\ref{ass.prob_good_sample}, so the desired conclusions follow.
\eproof

Let us close by observing that the proof of Theorem~\ref{th.appendix} does not require a measurable selection of the type in Assumption~\ref{ass.prob_good_sample}(a)--(b), since it merely employs the countable family of values $\{\beta_\phi/i\}_{i \in \N{}}$.

\section{Nonsmooth Rosenbrock Function}\label{app.rosenbrock}

In this appendix we report an experiment with a two-variable problem for which the objective, its generalized gradients, and its minimizer are all available in closed form, so that the behavior of Algorithm~\ref{alg.gs} can be compared directly with what our analysis predicts.  We considered a nonsmooth variant of the classical Rosenbrock function as proposed by Nesterov; see \cite{GurbOver2012}, namely, the function
\bequation\label{eq.rosenbrock}
  f(x,y) = (1 - x)^2 + 10|y - 2x^2 + 1|,
\eequation
which has the unique minimizer $(1,1)$ with $f(1,1) = 0$. We introduced artificial noise such that $\fapprox(x,y) = f(x,y) + \xi_f(x,y)$ and $\gapprox(x,y) = \nabla f(x,y) + \xi_g(x,y)$, where $\xi_f$ and $\xi_g$ are functions that are defined through a pseudo-random number generator whose seed is determined by the point at which the functions are evaluated: for each $(x,y)$, the value $\xi_f(x,y)$ is generated from a uniform distribution over the interval $[-\epsilon_f,\epsilon_f]$ and the value $\xi_g(x,y)$ is generated from a uniform distribution over the Euclidean ball with radius $\epsilon_g$ centered at the origin.  Hence, any two evaluations at the same point yield the same values, and the noise is deterministic and bounded as required by Assumptions~\ref{ass.f} and~\ref{ass.g}. Indeed, since $\|\gapprox(x,y) - \nabla f(x,y)\|_2 \leq \epsilon_g$ wherever $f$ is differentiable, one finds that Assumption~\ref{ass.g} holds with $\epsilon_x = 0$.  (We confirmed that the algorithm never requested a generalized gradient at a point of nondifferentiability of $f$.) We set $\epsilon_g=\sqrt{\epsilon_f}$ based on the common assumption that the gradient error is proportional to the square root of the error in the function values, as in finite-difference derivative approximations. We conducted tests across four noise levels, namely, $\epsilon_f \in \{10^{-3},10^{-4},10^{-5},10^{-6}\}$, in each case starting from the initial point $(-1.2,1)$ and choosing $\epsilonls = 2.1\epsilon_f$ so that $\epsilonls > 2\epsilon_f$, as required by Algorithm~\ref{alg.gs}.  Since the generalized gradients of $f$ are large relative to the distances over which they change (e.g., $\|\nabla f(-1.2,1)\|_2 \approx 53$), we chose $\nu = 10$ in light of the trade-off discussed in \S\ref{sec.implementation}; with the other parameters as stated there, this gives $\zeta = 2.37\epsilon_f^{1/2}$ in \eqref{eq.stuffff}, $\epsilon_1 = 23.7\epsilon_f^{1/2}$, and $\nu\zeta = 23.7\epsilon_f^{1/2} \geq \phi\epsilon_g$ for any $\phi \in (2,23]$.

\begin{figure}[htbp]
    \centering
    \includegraphics[width=0.6\textwidth]{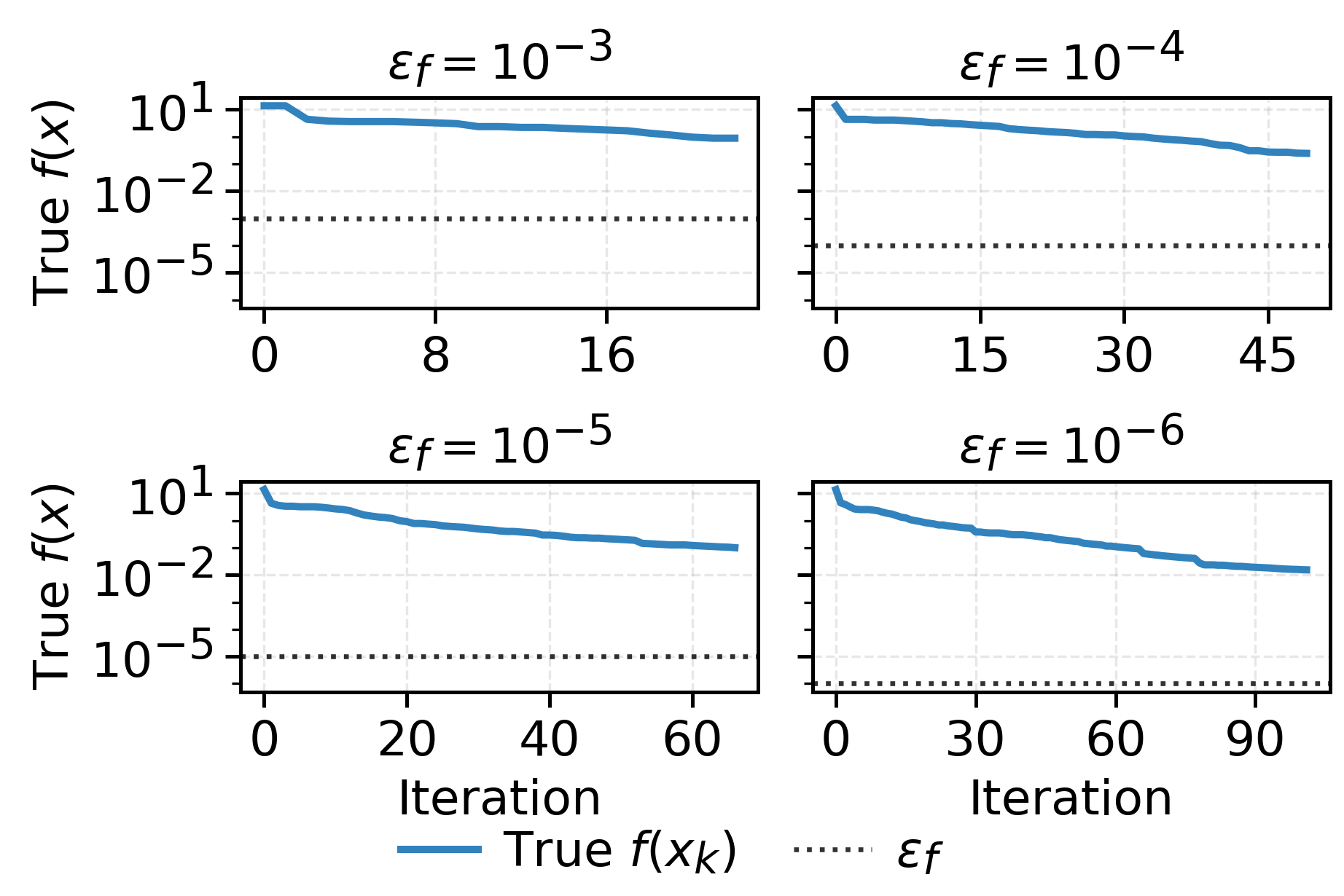}
    \caption{Nonsmooth Rosenbrock: True objective function values over $k$ obtained by NonOpt across the four noise levels $\epsilon_f \in \{10^{-3},10^{-4},10^{-5},10^{-6}\}$.  The dotted line in each plot marks $\epsilon_f$.}
    \label{fig:Rosen_2x2_nonsmooth}
\end{figure}

\begin{figure}[htbp]
    \centering
    \includegraphics[width=0.6\textwidth]{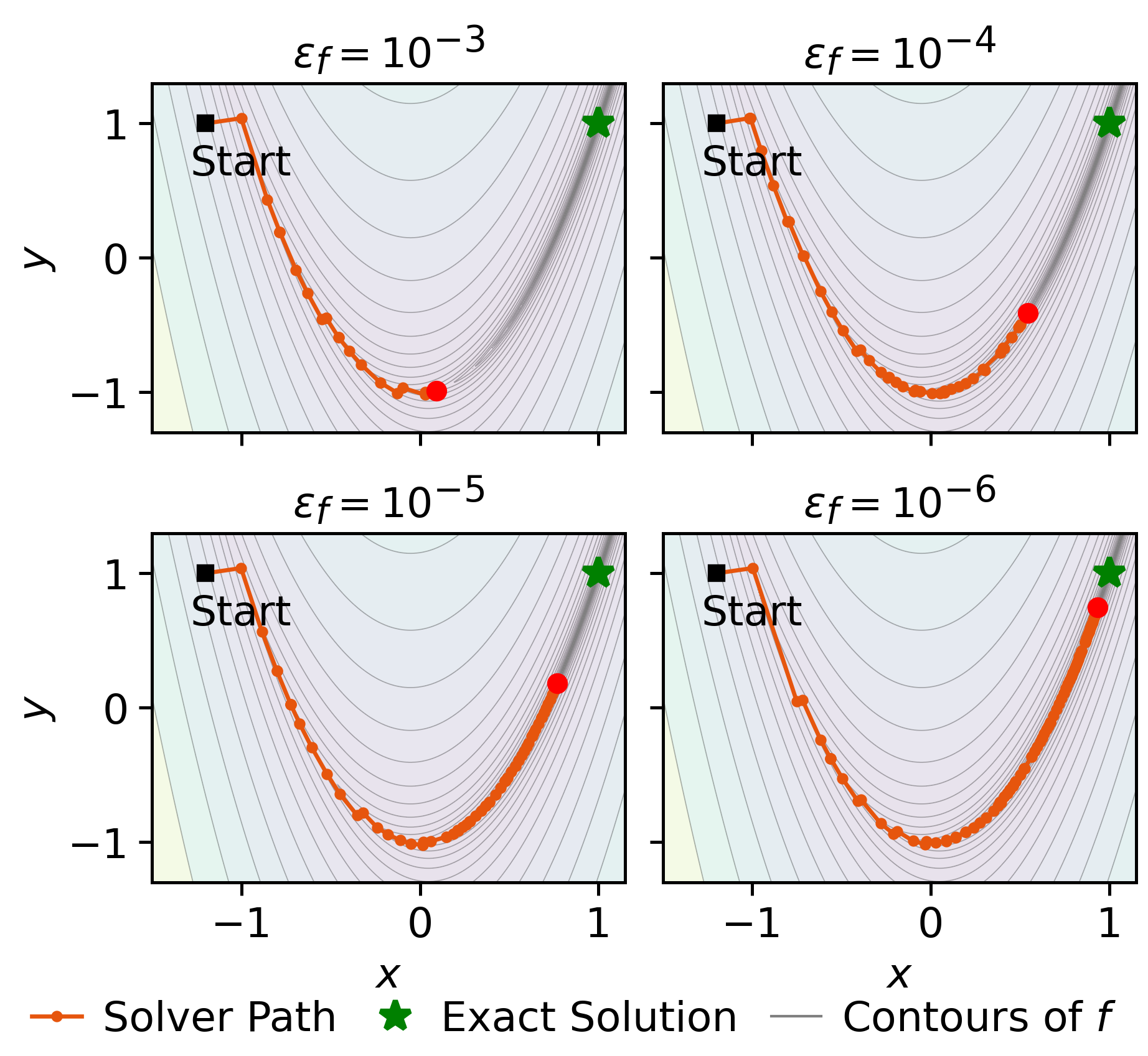}
    \caption{Nonsmooth Rosenbrock: Iterate trajectories obtained by NonOpt across the four noise levels $\epsilon_f \in \{10^{-3},10^{-4},10^{-5},10^{-6}\}$.  In each plot, the black square marks the initial point, the red dot marks the final iterate, and the green star marks the minimizer of $f$.}
    \label{fig:x_trajectory}
\end{figure}

Table~\ref{tab.rosenbrock} summarizes the runs, Figure~\ref{fig:Rosen_2x2_nonsmooth} shows the behavior of $\{f(x_k)\}$, and Figure~\ref{fig:x_trajectory} shows the corresponding iterate trajectories.  As Theorem~\ref{th.main-result} predicts, every run terminated finitely through the test in line~\ref{step.terminate} of Algorithm~\ref{alg.gs}.  Consistent with Lemma~\ref{lem.well-posed}(d), each of the $199$ steps with a positive step size that were taken over the four runs decreased $\fapprox$ by more than $\epsilonls$ (in fact, by at least $18\epsilonls$) and decreased the true objective $f$, so the true objective values in Figure~\ref{fig:Rosen_2x2_nonsmooth} decrease monotonically.  In all cases, the iterates reach the steep curved valley that contains the global minimizer, and, as the noise level decreases, they travel farther along the valley toward the minimizer before the algorithm terminates. (The first components of the final iterates were $0.093$, $0.54$, $0.77$, and $0.93$, respectively.)  As for stationarity, each run terminated with $\epsilon_k = \zeta$ and with a value of $\|\gapprox_k\|_2$ that is within the bound of Theorem~\ref{th.main-result}(b).  Moreover, in each case, the convex hull of the gradients of $f$ at the final iterate $x_k$ and at $1{,}000$ points sampled uniformly from $\ball{2}{2\epsilon_k}{x_k}$ contained the origin (to within $10^{-8}$), so $\dist_{\gengradient_{2\epsilon_k} f(x_k)}(0)$ is essentially zero, in agreement with the theorem.  The final objective values are nevertheless between $870$ and $1.6 \times 10^{4}$ times $\epsilon_f$.  This is a consequence of the tolerances that can be guaranteed for this problem rather than of the implementation: for any $x \in (\tfrac{1}{41},1)$, the point $(x,2x^2-1)$ on the valley floor is $\epsilon$-stationary in the sense of Goldstein for all $\epsilon \geq (1-x)/20$, since $\thalf \nabla f(x-t,2x^2-1) + \thalf \nabla f(x+t,2x^2-1) = 0$ for $t = (1-x)/20$, and $f(x,2x^2-1) = (1-x)^2$; hence, stationarity with respect to a sampling radius that is proportional to $\epsilon_f^{1/2}$ is compatible with objective values that are proportional to $\epsilon_f$ with a large constant.

\begin{table}[htbp]
  \footnotesize
  \caption{Nonsmooth Rosenbrock: results obtained by NonOpt for each noise level.  The columns are as in Table~\ref{tab.hinf} of \S\ref{sec.hinf}.}
  \label{tab.rosenbrock}
  \begin{center}
  \setlength{\tabcolsep}{5pt}
  \begin{tabular}{@{}cccccccccc@{}}
    \toprule
    $\epsilon_f$ & $\zeta$ & iterations & steps & abandoned & evaluations & $f$ & $\|\gapprox_k\|_2$ & $\nu\epsilon_k$ & bound \\
    \midrule
    $10^{-3}$ & $0.075$ & $23$ & $16$ & $5$ & $245$ & $0.88$ & $0.378$ & $0.748$ & $7.4$ \\
    $10^{-4}$ & $0.024$ & $50$ & $40$ & $8$ & $565$ & $0.24$ & $0.155$ & $0.236$ & $2.3$ \\
    $10^{-5}$ & $0.0075$ & $67$ & $56$ & $9$ & $710$ & $0.10$ & $0.0745$ & $0.0748$ & $0.74$ \\
    $10^{-6}$ & $0.0024$ & $102$ & $87$ & $13$ & $1{,}070$ & $0.015$ & $0.0215$ & $0.0236$ & $0.23$ \\
    \bottomrule
  \end{tabular}
  \end{center}
\end{table}

\section{Neural Network Training for Binary Classification}\label{app.nn}

In this appendix, to explore the behavior of our algorithm in a setting with stochastic noise---which is not covered by our analysis, but is common in practice---we consider the training of a neural network for binary classification when objective and generalized-gradient values are estimated by mini-batching and the noise bounds are not presumed to be known. We used the Cleveland database from the heart disease dataset~\cite{JanoSteiPfisDetr1989}; specifically, we used the records of $302$ of the $303$ patients in this database, each with $13$ features (normalized to have zero mean and unit variance) and a binary label that indicates the presence or absence of heart disease. The model was a feed-forward neural network with a single hidden layer of size 10 with ReLU activation and a single output logit, which is converted to a label probability via a sigmoid function, so that $n = 151$. The loss was binary cross-entropy. Due to the ReLU activation, the resulting objective is nonsmooth, while also nonconvex.

\begin{figure}[htbp]
  \centering
  \includegraphics[width=0.8\textwidth]{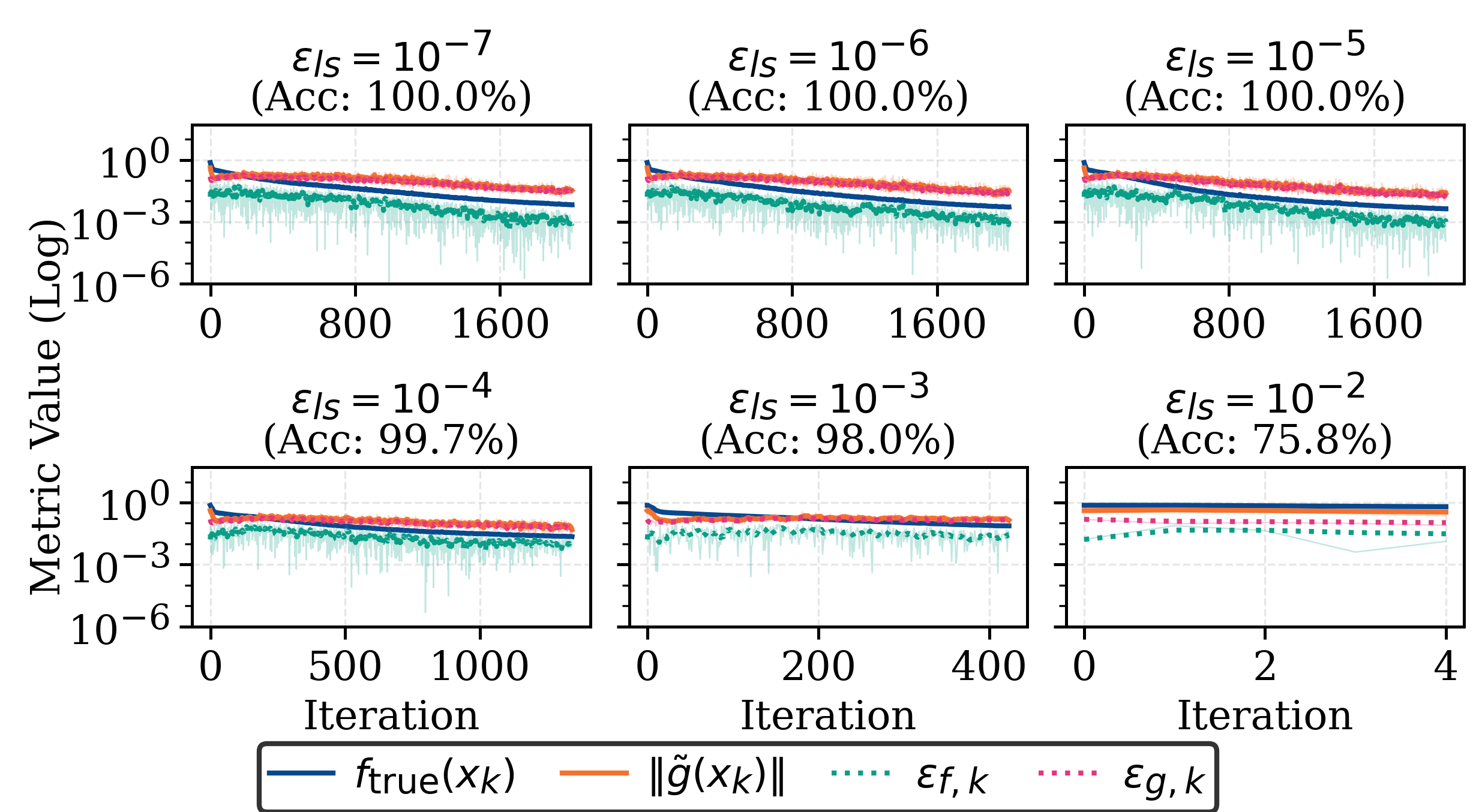}
  \caption{Neural network training with mini-batches of size 128: $f(x_k)$ (labeled $f_{\text{true}}$), $\|\gapprox(x_k)\|_2$, $\epsilon_{f,k}$, and $\epsilon_{g,k}$ over $k$ for different values of $\epsilonls$.  Faint traces are raw values; bold curves are moving averages with a window size of 8 iterations.  The title of each plot states the final training accuracy (``Acc.'').}
  \label{fig:fixed-mini-batch}
\end{figure}

Let $f$ and $g$ denote the objective and its (generalized) gradient as computed using the entire dataset, and let $\fapprox$ and $\gapprox$ denote the approximations of them that are computed using a mini-batch of 128 data points that is sampled anew in each iteration.  Specifically, all values of $\fapprox$ and $\gapprox$ that are employed in iteration $k$---including $\fapprox(x_k)$ and $\gapprox(x_k)$---are computed with the mini-batch of iteration $k$; in particular, the values of $\fapprox$ on the two sides of \eqref{eq.armijo} are computed with the same mini-batch.  We ran NonOpt as described in~\S\ref{sec.implementation} with $\nu = \gamma/3$ for $\epsilonls \in \{10^{-7},10^{-6},10^{-5},10^{-4},10^{-3},10^{-2}\}$ and an iteration limit of $2{,}000$.  Since no bounds $\epsilon_f$ and $\epsilon_g$ are presumed to be known, the requirement $\epsilonls > 2\epsilon_f$ cannot be verified and the observation in \S\ref{sec.implementation} about the term $\phi\epsilon_g$ does not apply; the remaining requirements of Algorithm~\ref{alg.gs} were imposed with $\epsilon_x = 0$, which gives $\zeta = 12.6\epsilonls^{1/2}$ in \eqref{eq.stuffff}, $\epsilon_1 = 126\epsilonls^{1/2}$, and $\nu\zeta = 2.11\epsilonls^{1/2}$.  For each $k$ we recorded $f(x_k)$, $\|\gapprox(x_k)\|_2$, and the errors $\epsilon_{f,k} := |\fapprox(x_k) - f(x_k)|$ and $\epsilon_{g,k} := \|\gapprox(x_k) - g(x_k)\|_2$. We emphasize that $f$ and $g$ were computed only in order to report the results; they were not employed by the algorithm in any way. The results are shown in Figure~\ref{fig:fixed-mini-batch} and Table~\ref{tab.heart}.

\begin{table}[htbp]
  \footnotesize
  \caption{Neural network training with mini-batches of size 128: results obtained by NonOpt for each value of $\epsilonls$.  Here, ``steps'' and ``abandoned'' are the numbers of iterations in which a positive step size was taken and in which the line search was abandoned by the test in line~\ref{line.threshold} of Algorithm~\ref{alg.gs}, respectively, and the last two columns state the value of $f$ and the training accuracy at the final iterate.  (At the initial point, these are $0.738$ and $59.3\%$.)}
  \label{tab.heart}
  \begin{center}
  \setlength{\tabcolsep}{5pt}
  \begin{tabular}{@{}ccccccccc@{}}
    \toprule
    $\epsilonls$ & $\zeta$ & $\nu\zeta$ & iterations & steps & abandoned & ended by & $f$ & accuracy \\
    \midrule
    $10^{-7}$ & $0.0040$ & $6.7 \times 10^{-4}$ & $1{,}999$ & $1{,}999$ & $0$ & limit & $0.00674$ & $100\%$ \\
    $10^{-6}$ & $0.013$ & $0.0021$ & $1{,}998$ & $1{,}997$ & $0$ & limit & $0.00523$ & $100\%$ \\
    $10^{-5}$ & $0.040$ & $0.0067$ & $1{,}998$ & $1{,}997$ & $0$ & limit & $0.00428$ & $100\%$ \\
    $10^{-4}$ & $0.13$ & $0.021$ & $1{,}343$ & $1{,}341$ & $0$ & line~\ref{step.terminate} & $0.0223$ & $99.7\%$ \\
    $10^{-3}$ & $0.40$ & $0.067$ & $424$ & $422$ & $0$ & line~\ref{step.terminate} & $0.0739$ & $98.0\%$ \\
    $10^{-2}$ & $1.3$ & $0.21$ & $5$ & $3$ & $0$ & line~\ref{step.terminate} & $0.503$ & $75.8\%$ \\
    \bottomrule
  \end{tabular}
  \end{center}
\end{table}

Over all runs, the errors $\epsilon_{f,k}$ had a median of $0.0069$ and a maximum of $0.11$, so the requirement $\epsilonls > 2\epsilon_f$ of Algorithm~\ref{alg.gs} would call for a value of $\epsilonls$ of more than $0.22$.  With such a value, however, one has $\nu\zeta \geq 0.99 > 0.68 \geq \|\gapprox(x_1)\|_2$, in which case the algorithm terminates at the initial point; i.e., with mini-batches of this size, the noise is too large for our guarantees to be informative for this problem.  All of the values of $\epsilonls$ that we consider are smaller.  For $\epsilonls = 10^{-2}$, $\epsilonls = 10^{-3}$, and $\epsilonls = 10^{-4}$, the algorithm terminated by the test in line~\ref{step.terminate} after $5$, $424$, and $1{,}343$ iterations, having reduced $f$ from $0.738$ to $0.503$, $0.0739$, and $0.0223$ and reached training accuracies of $75.8\%$, $98.0\%$, and $99.7\%$, respectively.  For $\epsilonls = 10^{-5}$, $\epsilonls = 10^{-6}$, and $\epsilonls = 10^{-7}$, the algorithm had not terminated when the iteration limit was reached; the final values of $f$ were $0.00428$, $0.00523$, and $0.00674$ and the final training accuracies were $100\%$, $100\%$, and $100\%$, respectively.  In these runs, a step was taken in at least $99.9\%$ of the iterations.  That the algorithm makes progress with values of $\epsilonls$ that are far smaller than the errors $\epsilon_{f,k}$ is due to the fact that the two values of $\fapprox$ in \eqref{eq.armijo} are computed with the same mini-batch; thus, each iteration is an iteration of the algorithm for the objective that is defined by the mini-batch of that iteration, with respect to which there is no noise.  On the other hand, since the mini-batch changes from one iteration to the next, a step that satisfies \eqref{eq.armijo} need not decrease $f$; over all runs, $f$ decreased in $53.5\%$ of the $7{,}759$ steps, although the final value of $f$ is smaller than the initial value in each run in which a step was taken.  Overall, the smaller the value of $\epsilonls$, the smaller the final objective value in the runs that terminated, as expected since the sampling radius $\zeta$ and the tolerance $\nu\zeta$ with which the algorithm terminates are proportional to $\epsilonls^{1/2}$; the smallest final objective value, namely, $0.00428$ (with accuracy $100\%$), was obtained with $\epsilonls = 10^{-5}$.

\section{Further Details of the Numerical Experiments}\label{app.details}

In this appendix, we provide further details of the experiments in \S\ref{sec.numerical}.

\subsection{Imposing the Test for Reducing the Sampling Radius}\label{app.details.implementation}

Let us describe how aspect~(v) in \S\ref{sec.implementation} was imposed.  By default, NonOpt reduces the sampling radius---or terminates, if a reduction would violate $\epsilon_{k+1} \geq \max\{\zeta,\epsilon_x\}$---when $\|\gapprox_k\|_\infty \leq \nu\rho_k\epsilon_k$, where
\bequationNN
  \rho_k := \frac{\max\{1,\|\gapprox(x_1)\|_\infty,\|\gapprox(x_k)\|_\infty\}}{\max\{1,\|\gapprox(x_1)\|_\infty\}} \geq 1.
\eequationNN
It applies this test first with $\gapprox(x_k)$ in place of $\gapprox_k$, before any points are sampled.  That test is not the one in line~\ref{step.nu}: it uses $\|\cdot\|_\infty$ in place of $\|\cdot\|_2$, it carries the factor~$\rho_k$, and, when it is satisfied before sampling, no points are sampled in that iteration.  We therefore disabled it, by setting its tolerance factor to zero, and imposed line~\ref{step.nu} through the routines that provide the values of $\fapprox$.  Those routines know $\gapprox_k$ as soon as the line search requests its first trial point, which is $x_k - \gapprox_k$ because the initial step size is one. They evaluate $\|\gapprox_k\|_2 \leq \max\{\nu\epsilon_k,\phi\epsilon_g\}$ there, and, when it holds, they decline every trial step, so that the iteration leaves the objective value unchanged.  NonOpt reduces the sampling radius, or terminates, after an iteration that did not decrease the objective value. This is its objective-similarity test, whose tolerance and whose limit we set to zero.  Since NonOpt evaluates that test at the beginning of the following iteration, each iteration of Algorithm~\ref{alg.gs} in which the test in line~\ref{step.nu} is satisfied is followed by one iteration of the software in which the radius is reduced (or the run terminates), in which no step is taken, and whose sample points enter no decision.  We do not count such an iteration as an iteration of Algorithm~\ref{alg.gs}.  In each run that terminated there were three such iterations of the software, namely, one after each of the two iterations in which the test in line~\ref{step.nu} was satisfied and one at the outset (the objective-similarity test of NonOpt is always satisfied in its first iteration, so we started it with the sampling radius $\theta^{-1}\epsilon_1$, which that iteration reduced to $\epsilon_1$); together they consumed three evaluations of $\fapprox$ and $3m$ evaluations of $\gapprox$, which are included in the numbers of evaluations that we report.  Conversely, in an iteration in which Algorithm~\ref{alg.gs} takes no step for a different reason---the test in line~\ref{line.threshold}, or one of the iterations just described---the routines report, at the step of negligible length that NonOpt then accepts, the largest floating-point number strictly below $\fapprox(x_k)$ rather than $\fapprox(x_k)$ itself, so that no such iteration reduces the sampling radius.  (The value at that point is not evaluated in either case.)  With this arrangement the sampling radius is reduced, and the algorithm terminates, exactly when $\|\gapprox_k\|_2 \leq \max\{\nu\epsilon_k,\phi\epsilon_g\}$ with $\gapprox_k$ computed after the $m$ points have been sampled, and never otherwise; in particular, Theorem~\ref{th.main-result}(b) applies to the final iterate as stated, without an additional factor.  In our experiments the parameters are such that $\nu\max\{\zeta,\epsilon_x\} \geq \phi\epsilon_g$ for some $\phi \in \R{}_{>2}$, so that the term $\phi\epsilon_g$ never decides the test, since $\epsilon_k \geq \max\{\zeta,\epsilon_x\}$ for all $k \in \N{}$.

\subsection{Fixed-Order $H_\infty$ Controller Design}\label{app.details.hinf}

Let us justify the bound that leads to \eqref{eq.hinf.epsf}.  For any unit vector $v$, the function $\omega \mapsto \|(A_K - i\omega I)v\|_2^2$ is a quadratic with unit leading coefficient.  If $v$ is a right singular vector corresponding to $\sigma_{\min}(A_K - i\omega^* I)$ for a minimizing frequency $\omega^*$, then this quadratic is at least $\sigma_{\min}(A_K - i\omega I)^2 \geq r(K)^2$ for all $\omega \in \R{}$ (since $A_K$ is real, $\sigma_{\min}(A_K - i\omega I)$ is even in $\omega$) and equals $r(K)^2$ at $\omega^*$, so it is minimized at $\omega^*$ and one finds that
\bequation\label{eq.hinf.quad}
  \sigma_{\min}(A_K - i\omega I)^2 \leq r(K)^2 + (\omega - \omega^*)^2\ \ \text{for all}\ \ \omega \in \R{}.
\eequation
As for the pruning rule with which each lattice minimum is evaluated, it is the following. Since $\sigma_{\min}(A_K - i\omega I)$ is $1$-Lipschitz in $\omega$, between two evaluated lattice points at distance $d$ with values $v_1$ and $v_2$ the values in between are at least $\thalf(v_1+v_2-d)$, so that stretch may be skipped once this exceeds the best value found.

Appendix~\ref{app.sampling} requires $\gapprox$ to be continuous on an open set of full measure.  Here $h$ is continuous and piecewise real-analytic in $K$, and $\gapprox$ is discontinuous only where the lattice point attaining the minimum changes or the smallest singular value at that lattice point has multiplicity greater than one. These sets lie in a closed set that is a locally finite union of zero sets of nonconstant real-analytic functions of $K$, hence of measure zero.  Consistently, at $2{,}000$ points sampled uniformly from $\ball{n}{10^{-3}}{K^*}$, a perturbation of size $10^{-9}$ changed $\gapprox$ by at most $2.9 \times 10^{-10}$ at every noise level, and the selected lattice point did not change at any of them.

\subsection{PDE-Constrained Optimization}\label{app.details.pde}

Let us state the targets that are used in \eqref{eq.pde.objective}.  With $b(\xi;a,\xi_0,s) := a \exp(-\|\xi - \xi_0\|_2^2/(2s^2))$ for $\xi \in \Omega$, the targets are the grid-point values of $u_d^{(1)} = b(\cdot;1,(0.30,0.35),0.15) + b(\cdot;0.35,(0.75,0.30),0.10)$ and $u_d^{(2)} = b(\cdot;1.1,(0.70,0.70),0.18) + b(\cdot;0.3,(0.35,0.80),0.09)$, so their main bumps lie in opposite corners of $\Omega$ and the two branches have well-separated basins.

A bound on the norms of the generalized gradients of the objective in \eqref{eq.pde.objective} is computable. Every element of $\gengradient f(x)$ can be expressed as $\sigma B^T J(u)^{-1} q + \lambda v$ for some $q$ with $\|q\|_2 \leq 1 + \mu$ and some $v \in [-1,1]^n$ (where the bound for $q$ is attained when the peak state value and the largest misfit occur at the same grid point with the state below the target), and since $J(u) \succeq A \succ 0$ one has $\|J(u)^{-1}\|_2 \leq \|A^{-1}\|_2$ for all $u$, so
\bequation\label{eq.pde.lipschitz}
  \|g\|_2 \leq \sigma \|B\|_2 \|A^{-1}\|_2 (1 + \mu) + \lambda \sqrt{n}\ \ \text{for all}\ \ g \in \gengradient f(x)\ \ \text{for all}\ \ x \in \R{n}.
\eequation
For the parameter values here, the right-hand side of \eqref{eq.pde.lipschitz} is approximately $5.15$.  (This bound is far from tight in the region that is relevant for our experiments. At the points at which generalized gradient approximations were computed in the runs reported in \S\ref{sec.pde}, their norms were between $0.065$ and $0.25$.)

\begin{figure}[ht]
  \centering
  \includegraphics[width=0.9\textwidth]{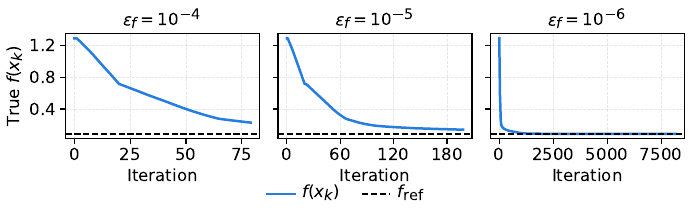}
  \caption{PDE-constrained problem: true objective value $f(x_k)$ over $k$ obtained by NonOpt across the three noise levels $\epsilon_f \in \{10^{-4},10^{-5},10^{-6}\}$, for the first of the five sampling seeds.  The dashed line in each panel marks the reference value $f_{\text{ref}} = 0.08334$.}
  \label{fig.pde.objective}
\end{figure}

The reference value $f_{\text{ref}}$ that is employed in \S\ref{sec.pde} is the smallest true objective value over the iterates of all of the runs that we performed.  These include a reference run, which used the same setup but with $\epsilon_f = 10^{-7}$, was started from the final iterate of the run with $\epsilon_f = 10^{-6}$, and was stopped at a time limit of $4$ hours.  The value $f_{\text{ref}}$ was attained at the final iterate of a separate run with the same setup except that $\nu = \gamma/9$ with $\epsilon_f = 10^{-6}$ (of all of the values of $\nu$ that we tested, this one gave the lowest objective value), which was stopped at its time limit of $2$ hours.  The runs reported in Table~\ref{tab.pde} required between $1.0 \times 10^{5}$ and $3.3 \times 10^{7}$ sweeps of the state iteration and between $3.8 \times 10^{5}$ and $5.1 \times 10^{7}$ CG iterations.  Figure~\ref{fig.pde.objective} shows the true objective values over the iterations of the runs for the first sampling seed.

\begin{figure}[t]
\centering
\includegraphics[width=0.8\textwidth]{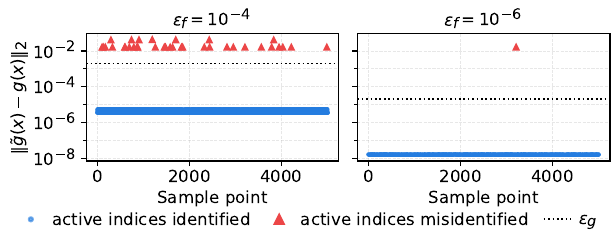}
\caption{PDE-constrained problem: generalized gradient error $\|\gapprox(x) - g(x)\|_2$ at $5{,}000$ points $x$ sampled uniformly from $\ball{n}{10^{-3}}{x_{\text{end}}}$ for $\epsilon_f = 10^{-4}$ (left) and $\epsilon_f = 10^{-6}$ (right).  Blue circles are points at which the inexact state selected the same active indices as the exact state; red triangles are those at which it did not.  The dotted line marks $\epsilon_g = 20\epsilon_f$.  The errors in the blue population decrease as the solver tolerances are tightened, whereas the errors in the red population do not; rather, the proportion of points in the red population decreases, illustrating the distinction between the $\epsilon_g$ and $\epsilon_x$ contributions in Assumption~\ref{ass.g}.}
\label{fig.pde.gradient}
\end{figure}

Let us provide further evidence for the distinction between the errors that are captured by $\epsilon_x$ and by $\epsilon_g$ in Assumption~\ref{ass.g}.  For each noise level, we evaluated $\gapprox$ at $5{,}000$ points sampled uniformly from a ball of radius $10^{-3}$ centered at the final iterate of the run with $\epsilon_f = 10^{-6}$ for the first sampling seed, which we denote by $x_{\text{end}}$.  This is a region in which $f$ is nonsmooth: $3$ distinct combinations of the active indices $(\ell^*,i^*,j^*)$ and sign $s$ occur at these points.  Figure~\ref{fig.pde.gradient} shows $\|\gapprox(x) - g(x)\|_2$ at these points, where $g(x)$ denotes the generalized gradient obtained by the same adjoint procedure with the state and adjoint equations solved to high accuracy, distinguishing the points at which the inexact state selected the same active indices and sign as the exact state (blue circles) from those at which it did not (red triangles).  At the former, the error was at most $0.078\epsilon_f$ in every case, which is well within $\epsilon_g$; i.e., the errors decreased as the solver tolerances were tightened.  At the latter, the error was between $0.016$ and $0.043$, which is of the order of the norms of the generalized gradients themselves (between $0.065$ and $0.072$ at these points), and did not decrease with the solver tolerances; what decreased is the proportion of such points, which was $0.70\%$, $0.12\%$, and $0.020\%$ for the three noise levels, respectively.  Had we attempted to bound the whole generalized gradient error by $\epsilon_g$ alone, we would have required $\epsilon_g \geq 0.016$ at every noise level, in which case the stationarity threshold in Theorem~\ref{th.main-result}, which depends additively on $\epsilon_g$, would not decrease however accurately the PDE is solved.

Appendix~\ref{app.sampling} requires $\gapprox$ to be continuous on an open set of full measure.  Here $\gapprox$ is discontinuous only where an iteration count of the state or adjoint solve changes, two terms of \eqref{eq.pde.objective} tie at $\utilde$, or some $x_j = 0$. These sets lie in a closed set, of measure zero unless two of the piecewise real-analytic quantities compared agree on an open set.  Consistently, at $3{,}000$ points sampled uniformly from $\ball{n}{10^{-3}}{x_{\text{end}}}$ and at every level, a perturbation of size $10^{-9}$ changed $\gapprox$ by more than $10^{-6}$ at none of them.